\documentclass[a4paper,11pt,reqno]{amsart}

\usepackage{amsmath, amsfonts, amssymb, amsthm, amscd}
\usepackage[foot]{amsaddr} % To put address at top page's foot
\usepackage[a4paper,scale={0.72,0.75},marginratio={1:1},footskip=7mm,headsep=10mm]{geometry}
\usepackage[T1]{fontenc}			
\usepackage[utf8]{inputenc}
\usepackage[italian, english]{babel}
\usepackage{csquotes}
\usepackage[square,numbers]{natbib}	
\usepackage{bbm}
\usepackage{bm}
\usepackage{mathtools}
\usepackage{xspace} % to get the spacing after macros right
\usepackage{booktabs}
\usepackage{perpage}
\usepackage{enumerate}
\usepackage{scalerel}
\usepackage{stackengine,wasysym}
\usepackage[pdftitle={Filtering of pure jump processes with partial observation},
pdfauthor={Elena Bandini, Giuseppina Guatteri, Gianmario Tessitore},%
]{hyperref}
\usepackage{verbatim} % Per commentare sezioni
\usepackage{mathrsfs}

\usepackage{stmaryrd}

\mathtoolsset{showonlyrefs}

\usepackage{hyperref}
\usepackage{dsfont}
\usepackage{color}

\newcommand{\R}{\mathbb{R}}

\newcommand{\N}{\mathbb{N}}

\renewcommand{\P}{\mathbb{P}}
\newcommand{\E}{\mathbb{E}}

\DeclarePairedDelimiterX{\abs}[1]{\lvert}{\rvert}{#1}
\DeclarePairedDelimiterX{\norm}[1]{\lVert}{\rVert}{#1}
\renewcommand{\epsilon}{\varepsilon}
\let\temp\theta
\let\theta\vartheta
\let\vartheta\temp

\newcommand{\phivarphi}
{
\let\temp\phi
\let\phi\varphi
\let\varphi\temp
}

\stackMath
\newcommand\reallywidehat[1]{%
\savestack{\tmpbox}{\stretchto{%
  \scaleto{%
    \scalerel*[\widthof{\ensuremath{#1}}]{\kern-.6pt\bigwedge\kern-.6pt}%
    {\rule[-\textheight/2]{1ex}{\textheight}}%WIDTH-LIMITED BIG WEDGE
  }{\textheight}%
}{0.5ex}}%
\stackon[1pt]{#1}{\tmpbox}%
}

\theoremstyle{plain}
\newtheorem{theorem}{Theorem}[section]
\newtheorem*{theorem*}{Theorem}
\newtheorem{lemma}[theorem]{Lemma}
\newtheorem*{lemma*}{Lemma}
\newtheorem{proposition}[theorem]{Proposition}
\newtheorem*{proposition*}{Proposition}
\newtheorem{corollary}[theorem]{Corollary}

\theoremstyle{definition}

\theoremstyle{remark}
\newtheorem{remark}{Remark}[section]
\newtheorem*{remark*}{Remark}

{\left\lbrace\begin{array}{@{}l@{}}}%
{\end{array}\right.}

\numberwithin{equation}{section}

\newcommand{\mail}[1]{\href{mailto:#1}{\normalfont\texttt{#1}}}

\makeatletter
\def\@setthanks{\vspace{-\baselineskip}\def\thanks##1{\@par##1\@addpunct.}\thankses}
\makeatother

\title[Singular limit of BSDEs and optimal control of two scale systems with jumps]
	{Singular limit of BSDEs and optimal control of two scale systems with jumps in infinite dimensional spaces}
\author[E.~Bandini]{Elena~Bandini}

\author[G.~Guatteri]{Giuseppina~Guatteri}

\author[G.~Tessitore]{Gianmario~Tessitore}

\date{}
\thanks{\noindent E.~Bandini
\\
University of Bologna, Department of Mathematics, piazza di Porta San Donato 5, 40126 Bologna (Italy).
\\
E-mail: \mail{elena.bandini7@unibo.it}.
\medskip
\\
G.~Guatteri
\\
Politecnico di Milano, Department of Mathematics, via Bonardi 9, 20133 Milan (Italy).
\\
E-mail: \mail{giuseppina.guatteri@polimi.it}.
\medskip
\\
G.~Tessitore
\\
University of Milano-Bicocca, Department of Mathematics and its Applications, via Roberto Cozzi 55, 20125 Milan (Italy).
\\
E-mail: \mail{gianmario.tessitore@unimib.it}.
\medskip
\\
This research was partially supported by the 2023 GNAMPA-INdAM  project \textit{Riduzione del modello in sistemi dinamici in dimensione infinita a due scale}.
\medskip
}

\begin{document}

\phivarphi
	
\begin{abstract} The paper is devoted to a stochastic optimal control problem for a two scale, infinite dimensional, stochastic system. The state of the system consists of “slow” and “fast” component and its evolution is driven by both continuous Wiener noises and discontinuous Poisson-type noises.  The presence of discontinuous noises is the main feature of the present work.  We use the theory of backward stochastic differential equations (BSDEs) to prove that,  as the speed of the fast component  diverges, the value function of the control problem converges to the  solution of a reduced forward backward system that, in turn, is related to a reduced stochastic optimal control problem. The results of this paper generalize to the case of discontinuous noise the ones  in \cite{guatess:amo2019}  and \cite{Swiech2021}.
\end{abstract}
\maketitle
	
\noindent \textbf{Keywords:} Discontinuous noise. Poisson random measures. Optimal control of two scales systems. Backward stochastic differential equations. Stochastic evolution equations in Hilbert spaces.

\smallskip

\noindent \textbf{AMS 2020:} 60H15, 60G51, 60G57, 93C70, 93E03.

\smallskip

\section{Introduction}

Stochastic dynamic models that evolve based on distinct and divergent time scales constitute a well-established area of research, encompassing both controlled and uncontrolled scenarios. This field has found extensive applications, as evidenced by classical results detailed in \cite{Kabanov_book}, and more recent applications, such as those explored in \cite{BardiKouhkouh} and \cite{FlandoliPappalettera}, particularly in the realms of neural networks and climate models.

Our specific focus lies in examining  two scales controlled systems featuring both Wiener and discontinuous (jump) noise:

\begin{align}\label{ControlledSDEintro}
\left\{
\begin{array}{llll}
d X^\varepsilon_t=   A X^\varepsilon_t dt  + b(X^\varepsilon_t, Q^\varepsilon_t, \alpha_t) dt +R d  W^1_t + \int_{H \setminus 0} w [N_1(dt\,dw) - {r(X^\varepsilon_t, Q^\varepsilon_t, \alpha_t,w)}\nu_1(dw) dt], \\
X_0 = x_0 \in H,\\
\varepsilon dQ^\varepsilon_t = [B Q_t^\varepsilon + F(X^\varepsilon_t, Q_t^\varepsilon)]dt + G \rho(\alpha_t)dt+ \sqrt \varepsilon G d  W^2_t\\ \displaystyle
\qquad \qquad +  \int_{K\setminus 0} w [\varepsilon N_2^\varepsilon(dt\,dw) - {\gamma(\alpha_t,w)}\nu_2(dw) dt],\\\displaystyle
Q_0^\varepsilon =q_0 \in K. 
\end{array}
\right.
\end{align}
In the above,  process $(\alpha)$ is the  control, $(X^\varepsilon)$ is the slowly evolving state variable and $(Q^\varepsilon)$ is the fast state variable. Notice that  $(\alpha)$, $(X^\varepsilon)$ and $(Q^\varepsilon)$  all take values in infinite dimensional spaces.
$(W^1_t)_{t \geq 0}$ and  $(W^2_t)_{t \geq 0}$) are independent {cylindrical} Wiener. 
$N_1(dt\,dw)$ and $N^\varepsilon_2 (dt \, dw)= N_2 (d({t}/{\varepsilon}) \, dw)$ are  Poisson random measures, {independent of each other and from $(W^1_t)_{t \geq 0}$ and  $(W^2_t)_{t \geq 0}$}, with  compensated random measures 
\begin{align*}
	\tilde N_1(dt\,dw)&=N_1(dt\,dw) - \nu_1(dw) dt,\\
 { \tilde{{N}}^\varepsilon_2(d(t/\epsilon)\,dw)}& {=N_2(dt\,dw) - \varepsilon^{-1}\nu_2(dw) dt},
\end{align*} while $r$ and $\gamma$ represent their controlled intensity. Finally, $A$ and $B$ are, possibly unbounded, linear operators,  we assume that $B+F$ is dissipative and that operator $R$ is invertible; in other words it is crucial to suppose that, if we forget the controlled terms, the fast variable has an ``ergodic'' behaviour while the slow variable is perturbed by non degenerate gaussian noise.

The above state equation is coupled with a standard cost functional as:
$$
J^{\varepsilon}(x_0, q_0, \alpha) = \E\left[\int_0^1 l(X^{\varepsilon}_t, Q_t^\varepsilon, \alpha_t)dt + h(X^{\varepsilon}_1)\right].
$$
The value function of the optimal control problem is then classically defined by:
\begin{equation}\label{valuefunc}
V^{\varepsilon}(x_0, q_0)=\inf_{\alpha} J^{\varepsilon}(x_0, q_0, \alpha).
\end{equation}
Our purpose is to characterize the limiting value function $\lim_{\varepsilon\to 0}V^{\varepsilon}(x_0, q_0)$.
%\textcolor{green}{In order to do that, we exploit the Backward Stochastic Differential Equation (BSDE) approach in this  framework to represent $V^\varepsilon(x_0, q_0)$ by means of a suitable forward-backward system. }
%
%\color{green}{In the diffusive framework,  stochastic optimal control problems of infinite dimensional processes have been
%extensively studied using the theory of BSDEs ; we mention in particular the  papers \cite{FuTess:AOP}, \cite{FuTesell} and the last chapter of
%the recent book \cite{fabbri:soc}, where a detailed discussion of the literature can be found. Starting from the seminal paper \cite{SuQue},   the BSDEs approach to  optimal control problems for discontinuous finite-dimensional processes has also been exploited, while  few contributions are available in the discontinuous infinite dimensional  framework, among which  \cite{bandiniconfortolacosso} in the case of the randomization method.
%\normalcolor

%\textcolor{magenta}{Several papers have been devoted to the above problem in the finite-dimensional case with continuous noise. For instance in \cite{Kabanov_book} and \cite{KabanovRunggaldier} results are obtained by  direct computations in specific situations while  in  \cite{AB3}, \cite{AB2}, \cite{AB1}, \cite{AC1}, \cite{AC2} the problem is studied through the convergence of the viscosity solutions to the Hamilton Jacobi Bellmen equation describing the evolution of the approximating value functions $V^{\varepsilon}$. }

Numerous papers have addressed the aforementioned issue within the finite-dimensional framework, specifically considering continuous noise. In works such as \cite{Kabanov_book} and \cite{KabanovRunggaldier}, results are derived through direct computations in particular scenarios. Meanwhile, in the investigations conducted in \cite{AB3}, \cite{AB2}, \cite{AB1}, \cite{AC1} and \cite{AC2}, the problem is approached by examining the convergence of viscosity solutions to the Hamilton-Jacobi-Bellman (HJB) equation, which describes the evolution of the approximating value functions $V^{\varepsilon}$.

% \textcolor{magenta}{As far as the infinite-dimensional case is concerned, to our best knowledge, the first results date back to \cite{guatess:amo2019} where the case of cylindrical Wiener process is treated assuming that the noise in the slow evolution is non-degenerate. In \cite{guatess:amo2019} a different approach is propose to overcome the technical difficulties embedded in the study of HJB equations in infinite dimensional spaces. Namely the approximating value functions $V^{\varepsilon}$ are represented by a forward-backward system of stochastic differential equations, see e.g. \cite{fabbri:soc}. The main technical point is then the proof that the solutions to the resulting singular sequence of farward-backward stochastic differential equations converge towards a reduced system living in a smaller space and involving the solution of an Ergodic Backward Stochastic Differential equation (for such a class of backward stochastic differential equations  see, for instance, \cite{FuhHuTess} or \cite{guatess:sicon2020}).}
Regarding the infinite-dimensional scenario, the earliest results, to the best of our knowledge, can be traced back to \cite{guatess:amo2019}. In this work, the cylindrical Wiener process case is addressed, assuming non-degeneracy of the noise in the slow evolution and introducing a new approach to overcome the technical challenges inherent in studying HJB equations in infinite-dimensional spaces. In particular, the approximating value functions $V^{\varepsilon}$ are represented through a forward-backward system of stochastic differential equations, as outlined in \cite[Chapter 6]{fabbri:soc}.

The critical technical aspect lies in proving that the solutions to the resulting singular sequence of forward-backward stochastic differential equations  converge towards a reduced system existing in a smaller space. This reduced system involves the solution of an ``ergodic''  BSDE, a class of equations discussed in works such as \cite{FuhHuTess} and \cite{guatess:sicon2020}.

%\textcolor{magenta}{ Then in \cite{guatess:sicon2022} the request on the non degeneracy of the noise in the slow equation is removed by a vanishing noise argument. The limit is finally represented as the value function of a reduced stochastic optimal control problem.}

Subsequently, in \cite{guatess:sicon2022}, the condition of non-degeneracy for the noise in the slow equation is eliminated through the application of a vanishing noise argument. The final outcome is expressed as the value function of a reduced stochastic optimal control problem.

Finally in \cite{Swiech2021} the author, assuming that the noise is gaussian and of trace class, successfully applies  viscosity solution techniques to the problem and  extends the range of treatable systems, including, for instance, the case of control dependent noise. 

%\textcolor{magenta}{Here the main novelty is the presence of discontinuous noise (both in the slow and in the fast equation). Following \cite{guatess:amo2019} we choose the forward backward stochastic system representation since it seems more adapted to the treatment  of non gaussian noises
%(since Poisson measures would induce integro-differential terms in the  HJB equation} 

In this context, the main innovation lies in the inclusion of discontinuous noise, present in both the slow and fast equations. Building on the approach introduced in \cite{guatess:amo2019}, we opt for the forward-backward stochastic system representation. This choice is motivated by its suitability for handling non-Gaussian noises, particularly because the inclusion of Poisson measures would introduce complicated integro-differential terms in the Hamilton-Jacobi-Bellman (HJB) equation,  see  \cite{SwiechZabczykUNIQUENESS},  \cite{SwiechZabczykEXISTENCE}.

%\textcolor{green}{E:    perch\'e non mettere dall'inizio il caso generale, ovvero con $\int_{\Xi \setminus 0} \Gamma_t(\omega, w) [N(\omega, dt\,dw) - \nu(dw) dt]$?}
We introduce the following forward-backward system to represent $V^\varepsilon(x_0, q_0)$ (see \cite{bandiniconfortolacosso} for the extension of the classical results to the discontinuous case):
\begin{align*}
\left\{
\begin{array}{llllll}
d X_t=   A X_t dt   +R d W^1_t + \int_{H \setminus 0} w [N_1(dt\,dw) - \nu_1(dw) dt], \\
X_0 = x_0 \in H,\\
\varepsilon dQ^\varepsilon_t = [B Q_t^\varepsilon + F(X_t, Q_t^\varepsilon)]dt + \sqrt \varepsilon G d W^2_t+ \int_{K\setminus 0} w [\varepsilon N_2^\varepsilon(dt\,dw) - \nu_2(dw) dt],\\
Q_0^\varepsilon =q_0 \in K,\\
-dY_t^{\varepsilon}  = \psi\Big (X_t, Q_t^{\varepsilon},Z^{\varepsilon}_t,  \frac{\Xi^{\varepsilon}_t}{\sqrt{\varepsilon}}, U^{\varepsilon}_t(\cdot), \frac{\Theta^{\varepsilon}_t(\cdot)}{\varepsilon}\Big ) dt -  Z^{\varepsilon}_t d W^1_t -  \Xi^{\varepsilon}_t d  W^2_t \\
	- \int_{H \setminus 0} U^{\varepsilon}_t(w) (N_1(dt \, dw) - \nu_1(dw) dt)-  \int_{K \setminus 0} \Theta^{\varepsilon}_t(w) \Big(N_2^\varepsilon(dt \, dw) - \frac{\nu_2(dw)}{\varepsilon} dt\Big), \\
	Y_1^{\varepsilon}=  h(X_1).
\end{array}
\right.
\end{align*}
The main technical result, see Theorem \ref{T:convYesp},  consists in the proof that the solution of the forward-backward system above (depending on $\varepsilon$ in a singular way)
converges towards the solution of the following ``reduced'' system \begin{align}\label{F-B_limit_system_2_inr
tro}
\left\{
\begin{array}{llll}
d X_t=   A X_t dt + R d W_t^1 + \int_{H \setminus 0} w [N_1(dt\,dw) - \nu_1(dw) dt], \\
X_0 = x_0,\\
 -d\bar Y_t =    \lambda(X_t,\bar Z_t,\bar U_t(\cdot)) dt 
	 -  \bar Z_t d  W^1_t -  \int_{H \setminus 0} \bar U_t(w) [N_1(dt \, dw) - \nu_1(dw) dt], \\
	\bar Y_1=h(X_1).
\end{array}
\right.
\end{align}
Here, $\lambda$ represents the value function of a suitable ergodic control problem. Similar to the approach in \cite{guatess:amo2019}, the proof relies on a discretization strategy for the slow variable and entails solving an ``ergodic'' BSDE, see \cite{FuhHuTess}. The introduction of Poisson random measures $N_1$ and $N_2$ introduces novel and non-trivial technical challenges.

%\normalcolor\rosso{[Qui ci vorrebbero dei riferimenti ragionati all'approccio BSDE ai problemi di controllo con salti (limitandosi al caso infinito dimensionale?)]}
%\color{green}{Stochastic optimal control problems of infinite dimensional processes have been
%extensively studied using the theory of BSDEs in the diffusive framework; we mention in particular the  papers \cite{FuTess:AOP}, \cite{FuTesell} and the last chapter of
%the recent book \cite{fabbri:soc}, where a detailed discussion of the literature can be found. Starting from the seminal paper \cite{SuQue},   the BSDEs approach to  optimal control problems for discontinuous finite-dimensional processes has also been exploited, while  few contributions are available in the discontinuous infinite dimensional  framework, among which  \cite{bandiniconfortolacosso} in the case of the randomization method.
%\\
%\textcolor{red}{E: per il caso di BSDEs che risolvono problemi di controllo in dimensione infinita con salti non sono riuscita a trovare altri riferimenti oltre il nostro lavoro \cite{bandiniconfortolacosso}. Forse c'è qualcosa di Oksendal con il principio del massimo, non so se e' il caso di citarlo... Sono anche in dubbio se questa sezione vada qui o prima, magari incorporata con due parole sulle BSDEs ergodiche e problemi di controllo ottimo. }
\normalcolor
%\textcolor{magenta}{Once the result on singular forward back ward system is proven we can complete the program representing the limit $\lim_{\varepsilon \to 0}V^{\varepsilon}(x_0, q_0)$ by the solution of a suitable reduced forward-backward system (see Theorem \ref{th:main_contr} and, eventually, as the value function of a suitable 'reduced' control problem  (see Theorem \ref{th:repr_as_control}) }.

After establishing the result concerning the singular forward-backward system, we can finalize the program by expressing the limit $\lim_{\varepsilon \to 0}V^{\varepsilon}(x_0, q_0)$ through the solution of an appropriate reduced forward-backward system (refer to Theorem \ref{th:main_contr}). Additionally, this limit can be further characterized as the value function of a suitable  ``reduced'' control problem (refer to Theorem \ref{th:repr_as_control}).

%\textcolor{magenta}{The paper is organized as follows: in Section 1 we set notations and general assumptions, in Section 2 we prove the main result on the convergence of singular forward-backward systems of stochastic equations, in Section 3 we apply the above convergence result to the original stochastic, two scales, control problem with jump-noise, finally, in Section 4, we interpret the main result in terms of a 'reduced' optimal control problem.}

The paper is structured as follows.  Section 2 establishes notations. In Section 3 we point out some results on the forward system. In Section 4 we present and prove the main result concerning the convergence of singular forward-backward systems of stochastic equations. Section 5 applies the aforementioned convergence result to the original stochastic two-scale control problem featuring jump noise. Finally, in Section 6, we provide an interpretation of the main result in the context of a ``reduced'' optimal control problem. We eventually provide some results on the Girsanov transformation for general semimartingales in the Appendix.

\section{Notation}\label{subsec:notation}
%In the next part we describe in detail the features  of our model and introduce the main assumptions and some specific notation.

%We recall that we are given a complete filtered probability space $(\Omega, \cF, \bbF \coloneqq (\cF_t)_{t \geq 0}, \bbP)$, that we assume rich enough to support all the probabilistic objects that we will encounter in the rest of this Section. We are also given a complete and separable metric space $(E, d_E)$ that we endow with the Borel $\sigma$-algebra $\cB(E)$.

For any topological space $\Lambda$,  $\mathcal B(\Lambda)$ denotes the corresponding Borel $\sigma$-algebra. 
Let $E$ be a Banach space norm $|\cdot|_E$ or simply $|\cdot|$ when no confusion is possible; for any other Banach space $F$, we denote by $L(E, F)$ the space of bounded linear operators from $E$ to $F$, endowed with the usual norm. If $F= \R$, the dual space  $L(E, \R)$ is denoted by $E^\ast$. We also use the shortwriting $L(E)$ for $L(E,E)$.

$\Xi, H, K$  are Hilbert spaces  with  scalar product denoted by $\langle \cdot, \cdot \rangle$. They  are assumed to be real and separable and the dual of Hilbert spaces will never be identified with the space itself. By $L_2(\Xi, H)$ (resp.  $L_2(\Xi, K)$) we will denote the space of Hilbert-Schmidt operators from $\Xi$ to $H$ (resp. to $K$). 
%\textcolor{green}{E: togliere - Moreover, $\mathcal G(K, H)$ will identify the space of all Gateaux differentiable mappings $\varphi$ from $K$ to $H$, such that the map $(k,v) \mapsto \nabla \varphi(k) v$  is continuous from $K \times K$ to $H$. }

We consider a complete probability space $(\Omega, \mathcal F, \P)$, and   a filtration $(\mathcal F_t)_{t \geq 0}$ satisfying the usual conditions. $T$ will be a fixed horizon, and 
$\mathcal P$ will denote the $\sigma$-algebra of predictable sets on $\Omega \times [0,\,T]$. 
 
{We consider two  cylindrical  $(\mathcal F_t)$-Wiener processes
$(W^1_t)_{t \geq 0}$ and  $(W^2_t)_{t \geq 0}$)  with values in $\Xi$.} We assume $(W^1_t)_{t \geq 0}$ and $(W^2_t)_{t \geq 0}$ to be independent.
$N_1(dt\,dw)$, $ {N}_2(dt\,dw)$ will be independent $(\mathcal F_t)$-Poisson random measures  on  $H$ and on $K$ respectively, with  compensated random measures 
\begin{align*}
	\tilde N_1(dt\,dw)&=N_1(dt\,dw) - \nu_1(dw) dt,\\
 \,\,\,\,\tilde{ {N}}_2(dt\,dw)& =N_2(dt\,dw) - \nu_2(dw) dt,
\end{align*}
where $\nu_1(dw)$ and  $\nu_2(dw)$ are  $\sigma$-finite measure on  $H$ and on $K$ respectively,  such that 
\begin{align*}
	\int_{H \setminus 0} |w|^2 \nu_1(dw) < \infty, \quad \int_{K \setminus 0} |w|^2 \nu_2(dw) < \infty.
\end{align*}
 $N_1(dt\,dw)$ and  $N_2(dt\,dw)$ are also independent of $(W^1_t)_{t \geq 0}$ and  $(W^2_t)_{t \geq 0}$.
Moreover for every $\varepsilon \in (0,1]$ we set
$$
N^\varepsilon_2 (dt \, dw)= N_2 (d({t}/{\varepsilon}) \, dw)
$$ 
that is a random measure 
  with compensator $\frac{1}{\varepsilon}\nu_2(dw) dt$.

We also introduce
\begin{align*}
L^2(\nu_1) &:= \{f: H \rightarrow \R  \,\,\textup{measurable}: \,\,\int_{H \setminus 0} |f(w)|^2\nu_1(dw) < \infty \},\\
L^2(\nu_2) &:= \{f: K \rightarrow \R   \,\,\textup{measurable}: \,\,\int_{K \setminus 0}  |f(w)|^2\nu_2(dw) < \infty \}, 
\end{align*}
and
\begin{align*}
	\mathbb{L}^2(\tilde N_1):=\{&\mathcal P \otimes \mathcal B(H) \textup{-measurable processes} \, U \,\,\textup{such that the norm }\\
	&||U||_{\mathbb{L}^2(\tilde N_1)} := \Big(\E\Big[\int_0^T\int_{H \setminus 0} |U_s(w)|^2 \nu_1(dw) ds\Big]\Big)^{1/2}\,\,\textup{is finite}\}, \\
	\mathbb{L}^2({\tilde{{N}}_2}):=\{&\mathcal P \otimes \mathcal B(K) \textup{-measurable processes}\, U \,\,\textup{such that the norm }\\
	&||U||_{\mathbb{L}^2({\tilde{{N}}_2})} := \Big(\E\Big[\int_0^T\int_{K \setminus 0} |U_s(w)|^2 \nu_2(dw) ds\Big]\Big)^{1/2}\,\,\textup{is finite}\}.
\end{align*}
	$\mathbb{L}^{2, \textup{loc}}(\tilde  {N}_1)$ (resp. $\mathbb{L}^{2, \textup{loc}}(\tilde N_2)$) denotes the set of processes defined on $\R_+$ such that their restriction to $[0,\,T]$ belongs to $\mathbb{L}^{2, \textup{loc}}(\tilde N_1)$ (resp. $\mathbb{L}^{2}(\tilde N_2)$).
 
 We also define the following spaces of equivalence classes of processes with values in a Hilbert space $V$:
%,  depending on a constant $p \geq 1$. 
	\begin{align*}
	S^2(\mathbb D([0,\,T]; V)):=\{& \mathcal{P}\textup{-measurable processes}\, X  \,\,\textup{with càdlàg paths in $V$, such that the   norm }\\
	&|X|_{S^2} := \Big(\E\Big[\sup_{ t \in [0,\,T]} |X_s|_V^2 ds\Big]\Big)^{1/2}\,\,\textup{is finite}\},\\
	 \mathbb{L}^2_V:= \{& \mathcal{P}\textup{-measurable processes}\, X  \,\,\textup{with values in $V$, such that the norm }\\
	& |X|_{\mathbb{L}^2_V} := \Big(\E\Big[\int_0^T |X_s|_V^2 ds\Big]\Big)^{1/2} \textup{is finite} \}.
\end{align*}
	$S^{2, \textup{loc}}(\mathbb D([0,\,+\infty); V))$ denotes the set of processes defined on $\R_+$ such that their restriction to $[0,\,T]$ belongs to $S^2(\mathbb D([0,\,T]; V)$.	
	$\mathbb{L}^{2, \textup{loc}}_V$ 
 	%(resp. $L^{2, \textup{loc}}(\mathbb D([0,\,T]; V))$)
 	 denotes the set of processes defined on $\R_+$ such that their restriction to $[0,\,T]$ belongs to $ \mathbb{L}^2_V$.

%	\verde{E:  possiamo toglierlo} 
%\begin{remark}
%\color{green}{According to \cite{PriolaZabczyk}, for every $U \in L^{2}(\tilde N_1)$ (resp. $\Theta \in L^{2}(\tilde N_2)$), we have 
%\begin{align*}
%	\E\Big[\Big|\int_0^t \int_{H \setminus 0} U_s(z) \tilde N_1(ds\,dz)\Big|^2\Big]&= ||U||^2_{\mathbb{L}^2(\tilde N_1)}\\
%	\Big(\textup{resp.} \,\,\E\Big[\Big|\int_0^t \int_{K \setminus 0} \Theta_s(w) \tilde N_2(ds\,dw)\Big|^2\Big]&= \frac{1}{\varepsilon}||\Theta||^2_{L^2(\tilde N_2)}\Big).
%\end{align*}}
%\end{remark}

\section{The forward system}
Given four operators   $A: D(A) \subseteq H \rightarrow H$, $B: D(B) \subseteq K \rightarrow K$, $R: \Xi \rightarrow H$,   $G: \Xi \rightarrow K$ and  $x_0 \in H$, $q_0 \in K$, we consider the following system of SDEs:
\begin{align}\label{SDE}
\left\{
\begin{array}{llll}
d X_t=   A X_t dt + R d W_t^1 + \int_{H \setminus 0} w [N_1(dt\,dw) - \nu_1(dw) dt], \\
X_0 = x_0,\\
\varepsilon dQ^\varepsilon_t = [B Q_t^\varepsilon + F(X_t, Q_t^\varepsilon)]dt + \sqrt \varepsilon G dW_t^2+  \int_{K\setminus 0}w [\varepsilon N_2^\varepsilon(dt\,dw) - \nu_2(dw) dt],\\
Q_0^\varepsilon =q_0, 
\end{array}
\right.
\end{align}
where the ``slow'' variable $X$ takes values in $H$ and the ``fast'' variable $Q^\varepsilon$ takes values in $K$.

We will consider the following assumptions on the coefficients of the system \eqref{SDE}.

\medskip 

\noindent {\sc{\textbf{(HAB)}}}
$A: D(A) \subseteq H \rightarrow H$ and $B: D(B) \subseteq K \rightarrow K$ are two linear unbounded  operators generating the $C_0$ semigroups $\{e^{tA}\}_{t\geq 0}$, $\{e^{tB}\}_{t \geq 0}$ over $H$ and $K$ respectively, such that 
\begin{align*}
	|e^{tA}|_{L(H)} \leq M_A e^{w_A t}, \quad t \geq 0, M_A >0,  w_A  \in \mathbb{R}, \\
	|e^{tB}|_{L(K)} \leq M_B e^{w_B t}, \quad t \geq 0, M_B >0 , w_B  \in \mathbb{R}.
\end{align*} 
Moreover, for $t>0$, $e^{tA}$ and $e^{tB}$ are Hilbert-Schmidt operators with 
$$
|e^{tA}|_{L_2(H,H)}+|e^{tB}|_{L_2(K,K)}\leq M t^{-\gamma t },
$$
 for all $t\in (0,1)$, suitable $\gamma \in (0,1/4)$ and $M>0$. 

\medskip 

\noindent {\sc{\textbf{(HRG)}}} $R \in L(\Xi, H)$, $G \in L(\Xi, K)$. 
\medskip 

\noindent {\sc{\textbf{(HF)}}} $F : H \times K \rightarrow K$ is bounded, and   there exists   $L_F>0$:  
\begin{align*}
	|F(x,y) - F(x', y')|_K \leq L_F (|x-x'|_H+|y-y'|_K), \quad x, x' \in H, y, y' \in K.
\end{align*}
%\verde{Moreover, for all $x \in H$, $F(x, \cdot) \in \mathcal G(K,K)$, SERVE? E: in effetti mi sembra di no}.

\noindent {\sc{\textbf{(HF+B)}}} $B+F$ is  dissipative: there exists $\mu >0$ such that 
\begin{align*}
	\langle B q + F(x,q) - [B q' + F(x,q')], q-q'\rangle \leq - \mu |q-q'|_K^2, \quad x\in H, q,q' \in D(B).
\end{align*}

Let $(\gamma_t)_{t \geq 0}$ be a  cylindrical Wiener process  with values in $\Xi$ and $p(\omega, dt\,dv)$  a random measure on  $K$ with compensator {$\kappa_t(\omega, v)\nu(dv) dt$, where 
$\nu(dw)$  is a $\sigma$-finite measure on $K$ and $\kappa$ is a random field on $K$   such that 
\begin{align}
\int_{K \setminus 0} |v|^2 |\kappa_s(\omega, v)|^2\nu(dv) < \infty \quad \textup{for every}\,\,(\omega,s) \in \Omega \times [0,T].\label{intknu}
\end{align}}
 We denote by $(\beta_t^B)_{t \geq 0}$ the stochastic convolution 
$$
\beta_t^B:= \int_0^t e^{(t-s)B}  \Big( G \, d \gamma_s + \int_{K \setminus 0} v [p(ds\,dv)- {\kappa_s(v)\nu(dv)} ds]\Big).
$$
We will assume the following:

\medskip

\noindent {\sc{\textbf{(H$\beta^B$)}}} $\sup_{t \geq 0} \E[|\beta_t^B|_K^2] < \infty$.
\begin{remark}
We notice that \textbf{(H$\beta^B$)}
is automatically satisfied when $B$ is dissipative and consequently $w_B<0$ in \textbf{(HAB).}
\end{remark}

\subsection{Well-posedness results}
The well-posedness of the slow equation in system \eqref{SDE} comes from the following theorem, that is a direct corollary of Theorem 9.29 in \cite{PeszatZabczyk}.
\begin{lemma}
Let assumptions  {\sc{\textbf{(HAB)}}}, {\sc{\textbf{(HRG)}}}, {\sc{\textbf{(HF)}}}, {\sc{\textbf{(HF+B)}}}, and  {\sc{\textbf{(H$\beta^B$)}}} hold. 	For any $T < \infty$,  $\tau \in [0,\,T]$ and any $\mathcal F_\tau$-measurable square integrable random variable $\tilde X_\tau \in H$,  
	 \begin{align*}
\left\{
\begin{array}{ll}
d X_t=   A X_t dt + R d W_t^1 + \int_{H \setminus 0} z [N_1(dt\,dz) - \nu_1(dz) dt], \quad t \in [\tau, T],\\
X_\tau = \tilde X_\tau
\end{array}
\right.
\end{align*}
has a unique (up to modification)  mild solution $(X_t^{\tau, \tilde X_\tau})_{t \in [\tau, T]}$ with a càdlàg version. More precisely, there exists a unique  process $(X_t)_{t \in [\tau, T]} \in S^2(\mathbb D([\tau,\,T]; H))$ such that 
\begin{align*}
	X_t &= e^{(t-\tau)A}\tilde X_\tau + \int_\tau^t e^{(s-\tau) A} R \, d W^1_s+ \int_\tau^t e^{(s-\tau) A}  \int_{H \setminus 0} z [N_1(ds\,dz) - \nu_1(dz) ds]
\end{align*}
is satisfied $\P$-a.s. Moreover, for every $T$ there exists $c_T <+\infty$ such that 
\begin{align*}
	\sup_{t \in [\tau, T]}\E[|X_t^{\tau, x}-X_t^{\tau, y}|_H^2] \leq c_T |x-y|_H^2, \quad x, y \in H,
\end{align*}
and, for all $p \geq 1$, there exists a positive constant $c_p$ depending only on $p$ and on the quantities introduced in the hypotheses,  such that  
\begin{equation}\label{crescitaLpX}
	\E\Big[\sup_{t \in [\tau, T]}|X_t^{\tau, x}|^p\Big]\leq c_p (1+ |x|^p), \quad x \in H. 
\end{equation}
\end{lemma}
Concerning  the fast equation in system \eqref{SDE},
we start by noticing that, if we make the change of time $s \mapsto s \varepsilon$, it reads 
\begin{align*}
\left\{
\begin{array}{ll}
\varepsilon d Q^{\varepsilon}_{\varepsilon t}=   [B Q^{\varepsilon}_{\varepsilon t} + F(X_t, Q^{\varepsilon}_{\varepsilon t})] \varepsilon dt  + \sqrt \varepsilon G dW_{\varepsilon t}^2+ \int_{K\setminus 0} w [\varepsilon N_2^\varepsilon(d(\varepsilon t)\,dw) - \varepsilon \nu_2(dw) dt], \quad t \in [0, T],\\
Q_0 = q_0,
\end{array}
\right.
\end{align*}
that in turn gives
\begin{align}\label{eqfast_changedtime}
\left\{
\begin{array}{ll}
d \hat Q_t=   [B  \hat Q_t + F(X_t, \hat Q_t)]  dt  +  G d\, \hat W^2_t +\int_{K\setminus 0} w [N_2(dt\,dw) - \nu_2(dw) dt], \quad t \in [0, T],\\
\hat Q_0 = q_0,
\end{array}
\right.
\end{align}
where we have set 
\begin{align*}
	\hat Q_t := Q^{\varepsilon}_{\varepsilon t}, \qquad
	\hat W^2_t := \frac{1}{\sqrt \varepsilon} dW_{\varepsilon t}^2.
	% N_2(dt\,dw)&:= N_2^\varepsilon(d(\varepsilon t)\,dw),
\end{align*}
% and $\hat L^2_t := \hat W^2_t + \int_{\Xi\setminus 0} z [N_2(dt\,dw) - \nu_2(dw) dt]$.
We have the following. 
\begin{lemma}\label{L:fasteq}
% Let assumptions  {\sc{\textbf{(HAB)}}}, {\sc{\textbf{(HRG)}}}, {\sc{\textbf{(HF)}}}, {\sc{\textbf{(HF+B)}}}, and  {\sc{\textbf{(H$\beta^B$)}}} hold.	
Let $(\Gamma_s)_{s \geq 0}$	be a given, $H$-valued, predictable process in $\mathbb{L}^{2, \textup{loc}}_H$. {Let $(g_s)_{s\geq 0}$ be a given, $K$-valued, process with $g \in \mathbb{L}^{2, \textup{loc}}_K$.}
%  for some $p\geq 1$.}
Let $(\gamma_t)_{t \geq 0}$ be a {$Q$}-Wiener process  with values in $\Xi$ and $p(\omega, dt\,dv)$  a random measure on  $K$ with compensator {$\kappa_t(\omega, v)\nu(dv) dt$, where 
$\nu(dw)$  is a $\sigma$-finite measure on $K$ and $\kappa$ is a random field on $K$   such that \eqref{intknu} holds true. 
%\begin{align*}
%\int_{K \setminus 0} |v|^2 |\kappa_t(\omega, v)|^2\nu(dv) < \infty \quad \textup{for every}\,\,(\omega,t) \in \Omega \times [0,T].
%\end{align*}
}
%Let $L$ be an $\Xi$-valued Lévy process with decomposition 
% 	$$
% 	L_t = \gamma_t + \int_0^t \int_{\Xi \setminus 0} z [p(dt\,dv) - \lambda(dv)dt],
% 	$$
Let us consider the following equation:
\begin{align}\label{eqQthm}
\left\{
\begin{array}{ll}
d Q_t=   [B Q_t + F(\Gamma_t, Q_t)]dt + g_t dt +  Gd\gamma_t + \int_{K \setminus 0}v[p(dt\,dv) - {\kappa_t(v)\nu(dv)}dt], \quad t \in [0, T],\\
Q_0 = q_0.
\end{array}
\right.
\end{align}
%where $p(\omega, dt\,dz)$ is a  random measure with compensator $\color{green}\kappa_t(\omega, v)\nu(dz)dt$ and $\gamma_t$ is acylindrical Wiener process.
 Then equation \eqref{eqQthm} admits a unique mild solution $Q \in S^{2, \textup{loc}}(\mathbb D([0,\,\infty);K))$. Moreover, for all $T>0$, there exists $k < \infty$, independent of $T$,  such that
\begin{align}\label{supestQ}
	\sup_{t \in [0, T]}\E[|Q_t|_K^2] \leq k \Big(1+ |q_0|_K^2  + \sup_{t \in [0, T]}\E[|\Gamma_t|_H^2]
	{+\sup_{t \in [0, T]}\E[|\beta^B_t|_K^2]}
	+\sup_{t \in [0, T]}\E[|g_t|_K^2]\Big).
\end{align}
Finally, if $(\Gamma'_s)_{s \geq 0}$ is another  process in $ \mathbb{L}^{2, \textup{loc}}_H$, and $Q'$ is a the mild solution of equation 
\begin{align*}
\left\{
\begin{array}{ll}
d Q'_t=   [B Q'_t + F(\Gamma'_t, Q'_t)]dt + g_t dt + Gd\gamma_t + \int_{K \setminus 0} v [p(dt\,dv) -{\kappa_t(v)\nu(dv)dt}], \quad t \in [0, T],\\
Q'_0 = q_0,
\end{array}
\right.
\end{align*}

then, for all $T >0$, $\P$-a.s., 
\begin{align}
    \label{diffestQ}
|Q_T - Q'_T|\leq \kappa\int_0^T e^{-\mu (T-s)}|\Gamma_s- \Gamma_s'|\,ds
\end{align}
where $\mu$ is the dissipativity constant in {\sc{\textbf{(HF+B)}}}, and $\kappa$ does not depend on $T$.
\end{lemma}
\proof
{The proof follows the same arguments of \cite{PeszatZabczyk} as far as existence and uniqueness  of the solution and estimate \eqref{supestQ} is concerned, while the last estimate follows as in \cite{DPZab2} (also see  \cite{guatess:amo2019}[Lemma 3.10]) noticing that
\begin{align*}
\left\{
\begin{array}{ll}
d (Q_t-Q'_t)=   [B(Q_t-Q'_t)+ F(\Gamma'_t, Q'_t)-F(\Gamma_t, Q_t)]dt , \quad t \in [0, T],\\ 
Q'_0 = q_0.
\end{array}
\right.
\end{align*}
% The well-posedness of \eqref{eqQthm} for every given process $\Gamma_s \in L^{2, \textup{loc}}([0,\,\infty);H)$ comes again from Theorem 9.29 in \cite{PeszatZabczyk}.
\endproof 
For fixed $\Gamma \equiv x \in H$, $q_0 \in K$ and $g \equiv 0$, equation \eqref{eqfast_changedtime} is a special case of equation \eqref{eqQthm} and Lemma \ref{L:fasteq} applies to it. We will denote by $\hat Q^{x, q_0}$ the mild solution to equation \eqref{eqfast_changedtime}.

%%%%%%%%%%%%%%%%%%%%%%%%%%%%%%%%%%%%%%%%%%%%%%%%%%%%%%%%

\section{Limit equation and convergence of singular BSDEs}
%\textcolor{red}{E: ho spostato qui l'ipotesi {\sc{\textbf{(H$h$)}}} che arrivava troppo presto, va uniformata con quello che c'è scritto sotto (mettere tutto in {\sc{\textbf{(H$h$)}}}?).}

We introduce a function $h:H \rightarrow \R$ satisfying the following assumption. 

\medskip 

\noindent 
 {\sc{\textbf{(H$h$)}}} 
$h:H \rightarrow \R $ is 
%  a progressively measurable and 
bounded, and 
 there exist positive constants   $L_h$ and $M$:   
\begin{align*}
	|h(x) - h(x')| \leq L_h |x-x'|_H,\quad |h(x)|  \leq M, \quad x, x' \in H, q, q' \in K.% $a \in U$,  
\end{align*}}

%The function $h$ is progressive measurable and bounded.}

\medskip 
We consider the following forward-backward system for $t \in [0,1]$:
\begin{align}\label{fbsystem_2}
\left\{
\begin{array}{llllll}
d X_t=   A X_t dt   +R d W^1_t + \int_{H\setminus 0} w [N_1(dt\,dw) - \nu_1(dw) dt], \\
X_0 = x_0 \in H,\\
\varepsilon dQ^\varepsilon_t = [B Q_t^\varepsilon + F(X_t, Q_t^\varepsilon)]dt + \sqrt \varepsilon G d W^2_t+  \int_{K\setminus 0} w [\varepsilon N_2^\varepsilon(dt\,dw) - \nu_2(dw) dt],\\
Q_0^\varepsilon =q_0 \in K,\\
-dY_t^{\varepsilon}  = \psi\Big (X_t, Q_t^{\varepsilon},Z^{\varepsilon}_t,  \frac{\Xi^{\varepsilon}_t}{\sqrt{\varepsilon}}, U^{\varepsilon}_t(\cdot), \frac{\Theta^{\varepsilon}_t(\cdot)}{\varepsilon}\Big ) dt -  Z^{\varepsilon}_t d W^1_t -  \Xi^{\varepsilon}_t d  W^2_t \\
	- \int_{H \setminus 0} U^{\varepsilon}_t(w) (N_1(dt \, dw) - \nu_1(dw) dt)-  \int_{K \setminus 0} \Theta^{\varepsilon}_t(w) \Big(N_2^\varepsilon(dt \, dw) - \frac{\nu_2(dw)}{\varepsilon} dt\Big), \\
	Y_1^{\varepsilon}=  h(X_1).
\end{array}
\right.
\end{align}
\begin{theorem}\label{T:welpos1BIS}
Assume  {\sc{\textbf{(HAB)}}}, {\sc{\textbf{(HRG)}}}, {\sc{\textbf{(HF)}}}, {\sc{\textbf{(HF+B)}}}, {\sc{\textbf{(H$\beta^B$)}}},   {\sc{\textbf{(H$\psi$)}}} and {\sc{\textbf{(H$h$)}}}. Then, for every $\varepsilon >0$,  there exists a unique 7-uple of processes $(X, Q^\varepsilon, Y^\varepsilon, Z^\varepsilon, \Xi^\varepsilon, U^\varepsilon(\cdot), \Theta^{\varepsilon}(\cdot))$ in $  S^2(\mathbb D([0,\,T]; H))\times S^2(\mathbb D([0,\,T]; K)) \times S^{2}(\mathbb D([0,\,+\infty); \R))\times \mathbb{L}^{2}_{\Xi^*}\times \mathbb{L}^{2}_{\Xi^*}\times  \mathbb{L}^{2}(\tilde{N}_1)\times  \mathbb{L}^{2}(\tilde{N}_2)$
	such that $\P$-a.s. the system  \eqref{fbsystem_2} is verified for all $t\in [0,1]$. 
\end{theorem}
\proof
 System \eqref{fbsystem_2} is a  decoupled forward-backward (that is, the forward equation can be solved independently, see for instance \cite{FuTess:AOP}).  The proof of existence and uniqueness  for the finite dimensional case can be found in \cite{confortola:filt} and can be extended easily to the infinite dimensional setting following \cite{bandiniconfortolacosso}.
\endproof

We aim at studying the limit behaviour of $Y^\varepsilon$ as $\varepsilon$ goes to $0$.

\subsection{The ergodic  parametrized BSDE}\label{S:erg} $ $

\noindent To formulate an equation for the limiting system we first have to introduce ergodic BSDEs in the context of discontinuous noise extending the results of \cite[section 4]{guatess:amo2019} and \cite[section 4]{FuhHuTess}. We remark that ergodic BSDEs with jumps in infinite dimension have also been studied in a different contest in  \cite{cohenErgodic}.

We assume that we are given a function $\psi: H \times K \times \Xi^\ast\times \Xi^\ast \times L^2(\nu_1) \times L^2(\nu_2) \rightarrow \R$, on which we ask the following.
 
 \medskip

\noindent {\sc{\textbf{(H$\psi$)}}}
\begin{itemize}
	\item[a)] $\psi$ measurable;
	\item[b)] $\sup_{x \in H, q \in K}\psi(x, q, 0,0,0,0) =:M_\psi < \infty$;
	\item[c)] there exist positive constants $L_x, L_q, L_z, L_\zeta,  L_u, L_\theta$ such that 
	 \begin{align*}
		&|\psi(x,q,z,\zeta,u, \theta)-\psi(x',q',z',\zeta',u', \theta')|\\
		&\leq L_z|z-z'|_{\Xi^\ast}+  L_\zeta|\zeta-\zeta'|_{\Xi^\ast}\\
		&+ L_u \Big(\int_{H \setminus 0}|u_s(w)- u'_s(w)|^2 \nu_1(dw)\Big)^\frac{1}{2} + L_\theta \Big(\int_{K \setminus 0}|\theta_s(w)- \theta'_s(w)|^2 \nu_2(dw)\Big)^\frac{1}{2}\\
		&+ L_x \Big(1 + |z|_{\Xi^\ast} + \Big(\int_{H \setminus 0}|u_s(w)|^2 \nu_1(dw)\Big)^\frac{1}{2}\Big)|x-x'|_H\\
		&+L_q \Big(1 + |z|_{\Xi^\ast} + \Big(\int_{H \setminus 0}|u_s(w)|^2 \nu_1(dw)\Big)^\frac{1}{2}\Big)|q-q'|_K;
	\end{align*}
	\item[d)] for every $x \in H$, $q \in K$, $z, \zeta \in \Xi^\ast$, $u, u' \in L^2(\nu_1)$, $\theta, \theta' \in L^2(\nu_2)$, there exist measurable functions $\gamma_1: H \rightarrow \R $, $\gamma_2: K \rightarrow \R $,  (that may depend on $x$, $q$, $z,\zeta$, $u$,  $u'$, $\theta, \theta'$) such that 
 \begin{align*}
		\psi(x,q,z,\zeta,u, \theta)-\psi(x,q,z,\zeta,u', \theta) \leq \int_{H \setminus 0} (u(w) -u'(w)) \gamma_1(w) \nu_1(dw),\\
		\psi(x,q,z,\zeta,u, \theta)-\psi(x,q,z,\zeta,u, \theta') \leq \int_{K \setminus 0} (\theta(w) -\theta'(w)) \gamma_2(w) \nu_2(dw),
	\end{align*}
	and satisfying 
	\begin{align*}
	C_1 (1 \wedge |w|_H) \leq \gamma_1(w) \leq C_2 (1 \wedge |w|_H), \quad C_1 \in (-1, 0], C_2 \geq 0, \\
	\bar C_1 (1 \wedge |w|_K) \leq \gamma_2(w) \leq \bar C_2 (1 \wedge |w|_K), \quad \bar C_1 \in (-1, 0], \bar C_2 \geq 0.
	\end{align*}
\end{itemize}
\begin{remark}
	Notice that the lipschitzianity of $f$ in the two latter components stated in \textbf{(H$f$)}-c) is a consequence of \textbf{(H$f$)}-d). 
 % Indeed, 
	% \begin{align*}
	% 	|\psi(x,q,z,\zeta,u, \theta)-\psi(x,q,z,\zeta,u, \theta')| &\leq C_2\int_{\Xi \setminus 0}(1+ ||w||_{\Xi}) |\theta(w) -\theta'(w)|  \nu_2(dw)\\
	% 	&\leq  L_\theta \Big(\int_{\Xi \setminus 0} |\theta(w) -\theta'(w)|^2  \nu_2(dw)\Big)^{\frac{1}{2}}. 
	% \end{align*}
	%with
	% $$
	% L_\theta=\Big(\int_{\Xi \setminus 0}(1+ ||w||_{\Xi})^2\nu_2(dw)\Big)^\frac{1}{2}
	% $$
	\end{remark}

\begin{remark}\label{Cohencond}
  %Following \cite{cohen:stochcalculus}[assumption 5] } 
  %\verde{[NON TROVO IL LAVORO SU MAHSCI] E: era il lavoro di Cohen che pero' non citerei qui, piuttosto diciamo che si verifica facilmente o mettiamo i passaggi }
  One can easily  deduce (see e.g. the proof of Lemma 19.3.8 in \cite{cohen:stochcalculus}) from 
 Hypothesis {\sc{\textbf{(H$\psi$)}}}-d) the  existence of two  
  functions 
  $\gamma_s^{1}$ and $\gamma_s^{2}$ (that may depend on $x$, $q$, $z,\zeta$, $u$,  $u'$ $\theta, \theta'$) satisfying 
\begin{align*}
	C_1 (1 \wedge |w|_H) \leq \gamma^{1}(w) \leq C_2 (1 \wedge |w|_H), \quad C_1 \in (-1, 0], C_2 \geq 0,\\
	\bar C_1 (1 \wedge |w|_K) \leq \gamma^{2}(w) \leq \bar C_2 (1 \wedge |w|_K), \quad \bar C_1 \in (-1, 0], \bar C_2 \geq 0,
\end{align*}
   %:=\gamma^{\hat {\mathcal{Q}}_s^{N},  \check \Xi^{N}_s, \hat \Theta^{\varepsilon}_s(\cdot), \check \Theta^{N}_s(\cdot)}$ 
     and such that 
  		\begin{align}
		&\psi(x,q,z,\zeta,u, \theta)-\psi(x,q,z,\zeta,u', \theta) = \int_{H  \setminus 0} \gamma^{1}(w)  (u(w) -u'(w)) \nu_1(dw)
	\end{align}
	and 
\begin{align}
		&\psi(x,q,z,\zeta,u, \theta)-\psi(x,q,z,\zeta,u(\cdot), \theta'(\cdot)) = \int_{K \setminus 0} \gamma^{2}(w)  (\theta(w) -\theta'(w)) \nu_2(dw).
	\end{align}
\end{remark}
We are interested in ergodic BSDEs of the form: $\P$-a.s., for all $x \in H$, $q \in K$, $z \in \Xi^\ast$, $u \in L^2(\nu_1)$,  
\begin{align}\label{ergodicBSDE}
	\hat Y_t &= \hat Y_T + \int_t^T [\psi(x, \hat Q_s^{x,q}, z, \hat \Xi_s, u(\cdot), \hat \Theta_s(\cdot)) - \lambda(x,z,u(\cdot))] ds \notag\\
	& - \int_t^T \hat \Xi_s d \hat W^2_s - \int_t^T \int_{\Xi \setminus  0} \hat\Theta_s(w) [N_2(ds \, dw) - \nu_2(dw) ds], \quad 0 \leq  t \leq T,
\end{align}
where $\hat Q_s^{x,q}$ is the unique solution to 
\begin{align}\label{eqhatQ}
\left\{
\begin{array}{ll}
d \hat Q_t=   [B  \hat Q_t + F(x, \hat Q_t)]  dt  + +  G d \hat W^2_t +\int_{K\setminus 0} w [N_2(dt\,dw) - \nu_2(dw) dt], \,\, t \in [0, T],\\
\hat Q_0 = q.
\end{array}
\right.
\end{align}
	\begin{remark}%\label{L:hatQ}
Under assumptions  {\sc{\textbf{(HAB)}}}, {\sc{\textbf{(HRG)}}}, {\sc{\textbf{(HF)}}}, {\sc{\textbf{(HF+B)}}}, and  {\sc{\textbf{(H$\beta^B$)}}}, 	the solution $\hat Q$  to 
	\eqref{eqhatQ} satisfies
		\begin{align}
			|\hat Q^{x,q}_t-\hat Q^{x,q'}_t|_K &\leq e^{-\mu t}|q-q'|_K, \quad t \in [0,\,T],\,\,q, q' \in K,\,x \in H,\label{esthatQ}
			\\
		%	\E\Big[\sup_{t \in [0,\,T]}|\hat Q_t^{x,q}|_K^p\Big] &\leq C_{p, T}(1+|q|_K^p), \quad q\in K,\,x \in H, \, p \geq 1,\label{estsupQ}
			%|\hat Q^{x,q}_t-\hat Q^{x,q'}_t|_K &\leq e^{-\mu t}|q-q'|_K, \quad t \in [0,\,T],\,\,q, q' \in K,\,x \in H,\label{esthatQ}
			%\\
			\E\Big[\sup_{t \in [0,\,T]}|\hat Q_t^{x,q}|_K^2\Big] &\leq C_{ T}(1+|q|_K^2), \quad q\in K,\,x \in H,\label{estsupQ}
		\end{align}
		for some constant $C_{T}$ only depending on $T$
(estimate \eqref{esthatQ} is a direct consequence of assumption {\sc{\textbf{(HF+B)}}}, while  the proof of \eqref{estsupQ} follows the same lines of the one of \eqref{supestQ}).
	\end{remark}
The theorem below is the analogue of the result in \cite[section 4]{guatess:amo2019}. We omit the proof since, on the one side it is similar to the one in \cite{guatess:amo2019} , on the other side the technical issues related to the presence of discontinuous noise will be treated, in all details, in the much more complicated proof of Theorem \ref{T:convYesp}.
\begin{theorem}\label{T:ergodic}
	Under assumptions   {\sc{\textbf{(HAB)}}}, {\sc{\textbf{(HRG)}}}, {\sc{\textbf{(HF)}}}, {\sc{\textbf{(HF+B)}}},  {\sc{\textbf{(H$\beta^B$)}}}, and {\sc{\textbf{(H$\psi$)}}}, there exist measurable functions 
	\begin{align*}
\bar v : H \times K \times \Xi^\ast \times L^2(\nu_1) \rightarrow \R,\\
\bar \theta : H \times K \times \Xi^\ast \times L^2(\nu_1) \rightarrow L^2(\nu_2),\\
\bar \zeta : H \times K \times \Xi^\ast \times L^2(\nu_1) \rightarrow \Xi^\ast,\\
\lambda : H  \times \Xi^\ast \times L^2(\nu_1) \rightarrow \R,
	\end{align*}
	such that 
	\begin{enumerate}
	\item for all fixed $x \in H$, $z \in \Xi^\ast$, $u \in L^2(\nu_1)$, $\bar v$ is Lipschitz with respect to $q$ and verifies, for some $c>0$, 	
	\begin{equation}\label{star1}
			|\bar v(x, q, z, u)|\leq c\Big(1+ |z| + \Big(\int_{ H\setminus 0}|u(\ell)|^2 \nu_1(d\ell) \Big)^\frac{1}{2}\Big) |q|_K;
	\end{equation}
\item if we set, for fixed  $x \in H$, $q \in K$,   $z \in \Xi^\ast$, $u \in L^2(\nu_1)$, 
\begin{align}\label{star2}
\left\{
\begin{array}{lll}
	\bar Y_t^{x, q, z, u}:= \bar v(x, \hat Q_t^{x, q}, z, u),\\
	\bar \Xi_t^{x, q, z, u}:= \bar \zeta(x, \hat Q_t^{x, q}, z, u),\\
	\bar \Theta_t^{x, q, z, u}(\cdot):= \bar \theta(x, \hat Q_t^{x, q}, z, u)(\cdot),
\end{array}
\right.
\end{align}
then $\bar \Xi^{x, q, z, u} \in \mathbb{L}^{2, \textup{loc}}$, $\bar \Theta^{x, q, z, u}(\cdot) \in \mathbb{L}^{2, \textup{loc}}(\tilde{N}_2)$, and the ergodic BSDE \eqref{ergodicBSDE} is satisfied by $(\bar Y^{x, q, z, u}, \bar \Xi^{x, q, z, u}, \bar\Theta^{x, q, z, u}(\cdot), \lambda(x, z,u))$.

Moreover, if in addition, we impose that there exists a constant $c>0$ such that
$
	|\bar Y_t^{x, q, z, u}|\leq c(1+|Q_t^{q_0}|),
$
then $\lambda(x, z, u)$ is uniquely determined. 

\item for fixed $x, x' \in H$, $z, z' \in \Xi^\ast$, $u, u' \in L^2(\nu_1)$, 
\begin{align}\label{star3}
	|\lambda(x,z,u)-\lambda(x',z',u')|&\leq L'_x\Big(1+ |z|_{\Xi^\ast} + \Big(\int_{H \setminus 0}|u(\ell)|^2 \nu_1(d\ell)\Big)^\frac{1}{2}\Big)|x-x'|_H\notag\\
	&+L'_z |z-z'|_{\Xi^\ast}+ L'_u\Big(\int_{H \setminus 0}|u(\ell)- u'(\ell)|^2 \nu_1(d\ell)\Big)^\frac{1}{2}
\end{align}
for some positive constants $L'_x$, $L'_z$, $L'_u$.
	\end{enumerate}
\end{theorem}

\subsection{Main Convergence Result}$ $

\noindent {We can now come back to our main object of investigation, namely the characterization of the limit, as $\varepsilon$ goes to $0$, of the sequence $Y_0^{\varepsilon} $ (where $Y^{\varepsilon}$ is the third component of the solution to  equation \eqref{fbsystem_2}.) }

 Our candidate limit is the component $\bar Y$ of the solution of the following forward-backward system on the finite horizon $[0,1]$ and on the reduced space $H$: 
\begin{align}\label{F-B_limit_system_2}
\left\{
\begin{array}{llll}
d X_t=   A X_t dt + R d W_t^1 + \int_{H \setminus 0} w [N_1(dt\,dw) - \nu_1(dw) dt], \\
X_0 = x_0,\\
 -d\bar Y_t =    \lambda(X_t,\bar Z_t,\bar U_t(\cdot)) dt 
	 -  \bar Z_t d  W^1_t -  \int_{H \setminus 0} \bar U_t(w) [N_1(dt \, dw) - \nu_1(dw) dt], \\
	\bar Y_1=h(X_1),
\end{array}
\right.
\end{align}
where $\lambda$ is the function  in Theorem \ref{T:ergodic} and $h$ is the function in {\sc{\textbf{(H$h$)}}}.

\begin{theorem}\label{T_wellposFB}
Let assumptions  {\sc{\textbf{(HAB)}}}, {\sc{\textbf{(HRG)}}}, {\sc{\textbf{(HF)}}}, {\sc{\textbf{(HF+B)}}}, {\sc{\textbf{(H$\beta^B$)}}},  {\sc{\textbf{(H$\psi$)}}} and  {\sc{\textbf{(H$h$)}}} hold true. Then, for every $\varepsilon >0$,  there exists a unique 4-uple of processes $( X,  \bar Y, \bar Z,  \bar U(\cdot))$ in $  S^2(\mathbb D([0,\,T]; H))\times S^{2}(\mathbb D([0,\,+\infty); \R))\times \mathbb{L}^{2}_{\Xi^*}\times  \mathbb{L}^{2}(\tilde{N}_1)$
	such that $\P$-a.s. the system  \eqref{F-B_limit_system_2} is verified for all $t\in [0,1]$.
Moreover,  there exist measurable functions 
	\begin{align*}
\bar v : H \rightarrow \R,\quad  
\bar \zeta : H  \rightarrow \Xi^\ast,\quad 
\bar \theta : H  \rightarrow L^2(\nu_1),
	\end{align*}
	such that
\begin{align*}
	\bar Y_t= \bar v(X_t),\quad 
	\bar Z_t= \bar \zeta(X_t),\quad 
	\bar U_t(\cdot)= \bar \theta(X_t)(\cdot).
\end{align*}
\end{theorem}
%\textcolor{red}{E: forse toglierei la proof, perché rimandiamo ad una proof prima che in realtà non facciamo.}
\proof 
 See the proof of Theorem \ref{T:welpos1BIS}.
\endproof

{The following is the main technical result of the paper.  The proof has the same structure of the one of \cite{guatess:amo2019}[Theorem 5.4]. The involvement of discontinuous noises induces nontrivial technical difficulties. The main differences lie in the application of the Girsanov transform (see Appendix \ref{AppGirs}) and in the estimates of the residual process $R^{\varepsilon, N}$, see \eqref{restoR}.}
\begin{theorem}\label{T:convYesp}
Let assumptions  {\sc{\textbf{(HAB)}}}, {\sc{\textbf{(HRG)}}}, {\sc{\textbf{(HF)}}}, {\sc{\textbf{(HF+B)}}}, {\sc{\textbf{(H$\beta^B$)}}},  {\sc{\textbf{(H$\psi$)}}} and  {\sc{\textbf{(H$h$)}}} hold true.
For every $\varepsilon >0$, let $(X, Q^\varepsilon, Y^\varepsilon, Z^\varepsilon, \Xi^\varepsilon, U^\varepsilon(\cdot))$ and $(X, \bar Y, \bar Z, \bar U(\cdot))$ be respectively the unique solution to the forward-backward systems \eqref{fbsystem_2} and \eqref{F-B_limit_system_2} obtained in Theorems \ref{T:welpos1BIS} and \ref{T_wellposFB}. Then 
$$
\lim_{\varepsilon \rightarrow 0}Y^\varepsilon_0= \bar Y_0.
$$
\end{theorem}
\proof
% We have 
% \begin{align*}
% 	Y^\varepsilon_0 -\bar Y_0&=\int_0^1 \Big(\psi\Big (X_t, Q_t^{\varepsilon},Z^{\varepsilon}_t,  \frac{\Xi^{\varepsilon}_t}{\sqrt{\varepsilon}}, U^{\varepsilon}_t(\cdot), \frac{\Theta^{\varepsilon}_t(\cdot)}{\varepsilon}\Big )-\lambda(X_t,\bar Z_t,\bar U_t(\cdot))\Big ) dt\\
% 	&-\int_0^1 (Z_t^\varepsilon - \bar Z_t) d W_t^1 -  \int_0^1\Xi^{\varepsilon}_t d  W^2_t\\
% 	&-\int_0^1\int_{\Xi \setminus 0} (U^{\varepsilon}_t(w)-\bar U_t(w)) (N_1(dt \, dw) - \nu_1(dw) dt)\\
% 	&-  \int_0^1\int_{\Xi \setminus 0} \Theta^{\varepsilon}_t(w) \Big(N_2^\varepsilon(dt \, dw) - \frac{\nu_2(dw)}{\varepsilon} dt\Big).
% \end{align*}

Adding and subtracting the term $\int_0^1 \psi\Big (X_t, Q_t^{\varepsilon},\bar Z_t,  \frac{\Xi^{\varepsilon}_t}{\sqrt{\varepsilon}}, \bar U_t(\cdot), \frac{\Theta^{\varepsilon}_t(\cdot)}{\varepsilon}\Big ) dt$ in the equation for $Y^\varepsilon$,  we have
\begin{align}\label{diff_to_est}
	Y^\varepsilon_0 -\bar Y_0&=\int_0^1 \Big(\psi\Big (X_t, Q_t^{\varepsilon},Z^{\varepsilon}_t,  \frac{\Xi^{\varepsilon}_t}{\sqrt{\varepsilon}}, U^{\varepsilon}_t(\cdot), \frac{\Theta^{\varepsilon}_t(\cdot)}{\varepsilon}\Big )-\psi\Big (X_t, Q_t^{\varepsilon},\bar Z_t,  \frac{\Xi^{\varepsilon}_t}{\sqrt{\varepsilon}}, \bar U_t(\cdot), \frac{\Theta^{\varepsilon}_t(\cdot)}{\varepsilon}\Big )\Big ) dt\notag\\
	&+\int_0^1 \Big(\psi\Big (X_t, Q_t^{\varepsilon},\bar Z_t,  \frac{\Xi^{\varepsilon}_t}{\sqrt{\varepsilon}}, \bar U_t(\cdot), \frac{\Theta^{\varepsilon}_t(\cdot)}{\varepsilon}\Big )-\lambda(X_t,\bar Z_t,\bar U_t(\cdot))\Big ) dt\notag\\
	&-\int_0^1 (Z_t^\varepsilon - \bar Z_t) d W_t^1 -  \int_0^1\Xi^{\varepsilon}_t d  W^2_t\notag\\
	&-\int_0^1\int_{H \setminus 0} (U^{\varepsilon}_t(w)-\bar U_t(w)) (N_1(dt \, dw) - \nu_1(dw) dt)\notag\\
	&-  \int_0^1\int_{K\setminus 0} \Theta^{\varepsilon}_t(w) \Big(N_2^\varepsilon(dt \, dw) - \frac{\nu_2(dw)}{\varepsilon} dt\Big).
\end{align}
Let $N \in \N$. We introduce a partition of the interval $[0,1]$ of the form $t_k = k 2^{-N}$, $k =0,1,..., 2^N$, and we define a triple of step processes $(X^N, \tilde Z^N, \tilde U^N(\cdot))$ as follows: 
\begin{align}
	X_t^N &= X_{t_k}, \quad t \in [t_k, t_{k+1}), \,\,k =0,..., 2^N-1,\label{XN}\\
	\tilde Z_t^N=0,\quad t \in [0, t_1),\,\,\textup{and}\,\, \tilde Z_t^N &= 2^N \int_{t_{k-1}}^{t_k}\bar Z_\ell d\ell, \quad t \in [t_k, t_{k+1}), \,\,k =1,..., 2^N-1,\label{ZN}\\
	\tilde U_t^N(\cdot)=0,\quad t \in [0, t_1),\,\,\textup{and}\,\, \tilde U_t^N(\cdot)&= 2^N \int_{t_{k-1}}^{t_k}\bar U_\ell^N(\cdot) d\ell, \quad t \in [t_k, t_{k+1}), \,\,k =1,..., 2^N-1.\label{UN}
\end{align}
By construction one has 
\begin{align}
	&\lim_{N \rightarrow \infty}\E \Big[\int_0^1 |\tilde Z^N_t - \bar Z_t|^2 dt\Big]=0, \label{convZN}\\
	&\lim_{N \rightarrow \infty}\E \Big[\int_0^1 \int_{H \setminus 0}|\tilde U^N_t(w) - \bar U_t(w)|^2 \nu_1(dw) dt\Big]=0.\label{convUN}
\end{align}
For $k=0,..., 2^N-1$ we consider the following class of forward SDEs of the type in  \eqref{eqfast_changedtime}:
\begin{align}
\left\{
\begin{array}{ll}\label{eqfast_changedtime_2}
d \hat {\mathcal{Q}}_t^{N,k}=   [B  \hat {\mathcal{Q}}^{N,k}_t + F(X_{t_k}, \hat {\mathcal{Q}}_t^{N,k})]  dt  + G  d\hat W^2_t + \int_{K\setminus 0} w [ N_2(dt\,dw) - \nu_2(dw) dt], \quad t \geq \frac{t_k}{\varepsilon},\\
\hat {\mathcal{Q}}^{N,k}_{\frac{t_k}{\varepsilon}} = \hat {\mathcal{Q}}^{N,k-1}_{\frac{t_k}{\varepsilon}},
\end{array}
\right.
\end{align}
%where $d\hat L^2_t = d\hat W^2_t + \int_{K\setminus 0} w [N_2(dt\,dw) - \nu_2(dw) dt]$ 
with
$\hat W^2_t = \frac{1}{\sqrt \varepsilon} dW_{\varepsilon t}^2$. 
At this point, let 
	\begin{align*}
\bar v : H \times K \times \Xi^\ast \times L^2(\nu_1) \rightarrow \R,\\
\bar \zeta : H \times K \times \Xi^\ast \times L^2(\nu_1) \rightarrow \Xi^\ast,\\
\bar \theta : H \times K \times \Xi^\ast \times L^2(\nu_1) \rightarrow L^2(\nu_2),
%\lambda : H  \times \Xi^\ast \times L^2(\nu_1) \rightarrow \R,
	\end{align*}
	be the functions introduced in Theorem \ref{T:ergodic}, and set, for $t \geq \frac{t_k}{\varepsilon}$,  
	\begin{align*}%\label{star2}
\left\{
\begin{array}{lll}
	\check  Y_t^{N,k}:= \bar v(X_{t_k}, \hat {\mathcal{Q}}_t^{N,k}, \tilde Z_{t_k}^N, \tilde U_{t_k}^N(\cdot)),\\
	\check \Xi_t^{N,k}:= \bar \zeta(X_{t_k}, \hat {\mathcal{Q}}_t^{N,k}, \tilde Z_{t_k}^N, \tilde U_{t_k}^N(\cdot)),\\
	\check \Theta_t^{N,k}(\cdot):= \bar \theta(X_{t_k}, \hat {\mathcal{Q}}_t^{N,k}, \tilde Z_{t_k}^N, \tilde U_{t_k}^N(\cdot))(\cdot).
\end{array}
\right.
\end{align*}
	By Theorem \ref{T:ergodic}, $(\check Y^{N,K}, \check \Xi^{N,k}, \check\Theta^{N,k}(\cdot), \lambda(X_{t_k}, \tilde Z_{t_k}^N,\tilde U_{t_k}^N(\cdot)))$
verifies, for all $k=0,..., 2^N-1$,
\begin{align}\label{ergodicBSDE_2}
	\check Y_{\frac{t_k}{\varepsilon}}^{N,k} &- \check Y^{N,k}_s - \int_{\frac{t_k}{\varepsilon}}^s [\psi(X_{t_k}, \hat {\mathcal{Q}}_s^{N,k}, \tilde Z_{t_k}^N, \check \Xi^{N,K}_s, \tilde U_{t_k}^N(\cdot), \check \Theta^{N,k}_s(\cdot)) - \lambda(X_{t_k},\tilde Z_{t_k}^N,\tilde U_{t_k}^N(\cdot))] ds \notag\\
	& + \int_{\frac{t_k}{\varepsilon}}^s \check \Xi^{N,k}_s d \hat W^2_s + \int_{\frac{t_k}{\varepsilon}}^s \int_{K \setminus 0} \check\Theta^{N,k}_s(w) [ N_2(ds \, dw) - \nu_2(dw) ds]=0, \quad s \geq \frac{t_k}{\varepsilon}. 
\end{align}
Moreover, there exists  
%	\begin{enumerate}
%	\item for all fixed $x \in H$, $z \in \Xi^\ast$, $u \in L^2(\nu_1)$, $\bar v$ is Lipschitz with respect to $q$ and verifies, 
 $c>0$ independent of $k$ and $N$, such that 	
	\begin{equation}\label{star1_2}
			|\check Y_s^{N,k}|\leq c\Big(1+ |\tilde Z_{t_k}^N| + \Big(\int_{H \setminus 0}|\tilde U_{t_k}^N(\cdot)|^2 \nu_1(d\ell) \Big)^\frac{1}{2}\Big) |\hat {\mathcal{Q}}_s^{N,k}|_K, \quad s \geq \frac{t_k}{\varepsilon}.
	\end{equation}
	Let us set, for $s \in [0,\frac{1}{\varepsilon}]$,  
	\begin{align}\label{mathcalQ^N}
		\hat {\mathcal{Q}}_s^N = \sum_{k=0}^{2^N-1}\hat {\mathcal{Q}}_s^{N,k} \,1_{\big[\frac{t_{k}}{\varepsilon},\frac{t_{k+1}}{\varepsilon}\big ]}(s),
	\end{align}
	and 
	\begin{align}\label{checkXiTheta}
	\check  {\Xi}_s^N = \sum_{k=0}^{2^N-1}\check  {\Xi}_s^{N,k} \,1_{\big[\frac{t_{k}}{\varepsilon},\frac{t_{k+1}}{\varepsilon}\big ]}(s),\quad \check  {\Theta}_s^N(\cdot) = \sum_{k=0}^{2^N-1}\check  {\Theta}_s^{N,k}(\cdot) \,1_{\big[\frac{t_{k}}{\varepsilon},\frac{t_{k+1}}{\varepsilon}\big ]}(s).
	\end{align}
	Taking $s = \frac{t_{k+1}}{\varepsilon}$ in \eqref{ergodicBSDE_2}, we get that, for all $k=0,..., 2^N-1$,
\begin{align}\label{ergodicBSDE_3}
	\check Y_{\frac{t_k}{\varepsilon}}^{N,k} &- \check Y^{N,k}_\frac{t_{k+1}}{\varepsilon} - \int_{\frac{t_k}{\varepsilon}}^\frac{t_{k+1}}{\varepsilon} \Big(\psi(X^N_{\varepsilon s}, \hat {\mathcal{Q}}_s^{N}, \tilde Z_{\varepsilon s}^N, \check \Xi^{N}_s, \tilde U_{\varepsilon s}^N(\cdot), \check \Theta^{N}_s(\cdot)) - \lambda(X^N_{\varepsilon s},\tilde Z_{\varepsilon s}^N,\tilde U_{\varepsilon s}^N(\cdot))\Big) ds \notag\\
	& + \int_{\frac{t_k}{\varepsilon}}^\frac{t_{k+1}}{\varepsilon} \check \Xi^{N}_s d \hat W^2_s + \int_{\frac{t_k}{\varepsilon}}^\frac{t_{k+1}}{\varepsilon} \int_{K \setminus 0} \check\Theta^{N}_s(w) [ N_2(ds \, dw) - \nu_2(dw) ds]=0,  
\end{align} 
	where we have 
set, for $s \in [0, \frac{1}{\varepsilon}]$, 
\begin{align*}
	\hat {\mathcal{Q}}_s^{N}&:= \sum_{k=0}^{2^N-1}\hat {\mathcal{Q}}_s^{N,k} \,1_{[\frac{t_k}{\varepsilon}, \frac{t_{k+1}}{\varepsilon})}(s),\\
	\check \Xi^{N}_s&:= \sum_{k=0}^{2^N-1} \check \Xi^{N,k}_s \,1_{[\frac{t_k}{\varepsilon}, \frac{t_{k+1}}{\varepsilon})}(s),\\
	\check \Theta^{N}_s&:= \sum_{k=0}^{2^N-1} \check \Theta^{N,k}_s \,1_{[\frac{t_k}{\varepsilon}, \frac{t_{k+1}}{\varepsilon})}(s), 
\end{align*}
and we have used that, for  $s \in [\frac{t_k}{\varepsilon}, \frac{t_{k+1}}{\varepsilon}]$, one has  $X_{t_k}= X^N_{\varepsilon s}$, $\tilde Z_{t_k}^N= \tilde Z_{\varepsilon s}^N$ and $\tilde U_{t_k}^N(\cdot)=\tilde U_{\varepsilon s}^N(\cdot)$ (recall \eqref{XN}-\eqref{ZN}-\eqref{UN}).

Let us now go back to \eqref{diff_to_est}. The term 
$$
\int_0^1 \Big(\psi\Big (X_t, Q_t^{\varepsilon},\bar Z_t,  \frac{\Xi^{\varepsilon}_t}{\sqrt{\varepsilon}}, \bar U_t(\cdot), \frac{\Theta^{\varepsilon}_t(\cdot)}{\varepsilon}\Big )-\lambda(X_t,\bar Z_t,\bar U_t(\cdot))\Big ) dt
$$
can be rewritten as  
\begin{align*}
	&\varepsilon\sum_{k=0}^{2^N-1}\int_{\frac{t_k}{\varepsilon}}^{\frac{t_{k+1}}{\varepsilon}} \Big(\psi (X_{\varepsilon s}, \hat Q_s^{\varepsilon},\bar Z_{\varepsilon s},  \hat \Xi_{s}^{\varepsilon}, \bar U_{\varepsilon s}(\cdot), \hat \Theta^{\varepsilon}_s(\cdot))-\lambda(X_{\varepsilon s},\bar Z_{\varepsilon s},\bar U_{\varepsilon s}(\cdot))\Big ) ds, 
\end{align*}
where we have denoted  $\hat Q_s^{\varepsilon} = Q^{\varepsilon}_{\varepsilon t}$, $\hat \Xi_s^{\varepsilon}=\frac{1}{\sqrt{\varepsilon}}\Xi^{\varepsilon}_{\varepsilon s}$ and $\hat \Theta_s^{\varepsilon}(\cdot)= \frac{1}{\varepsilon}\Theta^{\varepsilon}_{\varepsilon s}(\cdot)$.
Adding in the previous expression  the null terms in \eqref{ergodicBSDE_3}  for $k=0,..., 2^N-1$, we get 
\begin{align}\label{secondterm}
&\int_0^1 \Big(\psi\Big (X_t, Q_t^{\varepsilon},\bar Z_t,  \frac{\Xi^{\varepsilon}_t}{\sqrt{\varepsilon}}, \bar U_t(\cdot), \frac{\Theta^{\varepsilon}_t(\cdot)}{\varepsilon}\Big )-\lambda(X_t,\bar Z_t,\bar U_t(\cdot))\Big )\Big ) dt\notag\\
	&=\varepsilon\sum_{k=0}^{2^N-1}\Big  (\check Y_{\frac{t_k}{\varepsilon}}^{N,k} - \check Y^{N,k}_\frac{t_{k+1}}{\varepsilon}\Big )\notag\\
	&+\varepsilon\sum_{k=0}^{2^N-1}\int_{\frac{t_k}{\varepsilon}}^{\frac{t_{k+1}}{\varepsilon}} \Big(\psi (X_{\varepsilon s}, \hat Q_s^{\varepsilon},\bar Z_{\varepsilon s},  \hat \Xi^{\varepsilon}_{s}, \bar U_{\varepsilon s}(\cdot), \hat \Theta^{\varepsilon}_s(\cdot))-\psi(X^N_{\varepsilon s}, \hat {\mathcal{Q}}_s^{N}, \tilde Z_{\varepsilon s}^N, \check \Xi^{N}_s, \tilde U_{\varepsilon s}^N(\cdot), \check \Theta^{N}_s(\cdot))\Big ) ds\notag\\
	& - \varepsilon\sum_{k=0}^{2^N-1}\int_{\frac{t_k}{\varepsilon}}^\frac{t_{k+1}}{\varepsilon} \Big(\lambda(X_{\varepsilon s},\bar Z_{\varepsilon s},\bar U_{\varepsilon s}(\cdot))  - \lambda(X_{\varepsilon s}^N,\tilde Z_{\varepsilon s}^N,\tilde U_{\varepsilon s}^N(\cdot))\Big) ds \notag\\
	& + \varepsilon\sum_{k=0}^{2^N-1}\int_{\frac{t_k}{\varepsilon}}^\frac{t_{k+1}}{\varepsilon} \check \Xi^{N}_s d \hat W^2_s + \varepsilon\sum_{k=0}^{2^N-1}\int_{\frac{t_k}{\varepsilon}}^\frac{t_{k+1}}{\varepsilon} \int_{K \setminus 0} \check\Theta^{N}_s(w) (N_2(ds \, dw) - \nu_2(dw) ds).
\end{align}
On the other hand, the first term in the right-hand side of \eqref{diff_to_est} reads 
\begin{align}\label{firstterm}
&\int_0^1 \Big(\psi\Big (X_t, Q_t^{\varepsilon},Z^{\varepsilon}_t,  \frac{\Xi^{\varepsilon}_t}{\sqrt{\varepsilon}}, U^{\varepsilon}_t(\cdot), \frac{\Theta^{\varepsilon}_t(\cdot)}{\varepsilon}\Big )-\psi\Big (X_t, Q_t^{\varepsilon},\bar Z_t,  \frac{\Xi^{\varepsilon}_t}{\sqrt{\varepsilon}}, \bar U_t(\cdot), \frac{\Theta^{\varepsilon}_t(\cdot)}{\varepsilon}\Big )\Big ) dt	\\
&=\varepsilon\int_0^{\frac{1}{\varepsilon}} \Big(\psi(X_{\varepsilon s}, \hat Q_s^{\varepsilon},Z^{\varepsilon}_{\varepsilon s},  \hat \Xi^{\varepsilon}_{s}, U^{\varepsilon}_{\varepsilon s}(\cdot), \hat \Theta^{\varepsilon}_s(\cdot))-\psi(X_{\varepsilon s}, \hat Q_s^{\varepsilon},\bar Z_{\varepsilon s},  \hat \Xi^{\varepsilon}_s, \bar U_{\varepsilon s}(\cdot), \hat \Theta^{\varepsilon}_s(\cdot))\Big ) ds. \notag
\end{align}
Plugging \eqref{firstterm} and \eqref{secondterm} in \eqref{diff_to_est} we obtain
\begin{align}\label{diff_to_est_2}
	Y^\varepsilon_0 -\bar Y_0&=\varepsilon\sum_{k=0}^{2^N-1}\Big  (\check Y_{\frac{t_k}{\varepsilon}}^{N,k} - \check Y^{N,k}_\frac{t_{k+1}}{\varepsilon}\Big )\notag\\
	&+\varepsilon\int_0^{\frac{1}{\varepsilon}} \Big(\psi (X_{\varepsilon s}, \hat Q_s^{\varepsilon},Z^{\varepsilon}_{\varepsilon s},  \hat \Xi^{\varepsilon}_{s}, U^{\varepsilon}_{\varepsilon s}(\cdot), \hat \Theta^{\varepsilon}_s(\cdot) )-\psi (X_{\varepsilon s}, \hat Q_s^{\varepsilon},\bar Z_{\varepsilon s},  \hat \Xi^{\varepsilon}_s, \bar U_{\varepsilon s}(\cdot), \hat \Theta^{\varepsilon}_s(\cdot) )\Big ) ds\notag\\
	&+\varepsilon\int_0^{\frac{1}{\varepsilon}}\Big(\psi (X_{\varepsilon s}, \hat Q_s^{\varepsilon},\bar Z_{\varepsilon s},  \hat \Xi^{\varepsilon}_{s}, \bar U_{\varepsilon s}(\cdot), \hat \Theta^{\varepsilon}_s(\cdot) )-\psi(X^N_{\varepsilon s}, \hat {\mathcal{Q}}_s^{N}, \tilde Z_{\varepsilon s}^N, \check \Xi^{N}_s, \tilde U_{\varepsilon s}^N(\cdot), \check \Theta^{N}_s(\cdot))\Big ) ds\notag\\
	& - \varepsilon\int_0^{\frac{1}{\varepsilon}} \Big(\lambda(X_{\varepsilon s},\bar Z_{\varepsilon s},\bar U_{\varepsilon s}(\cdot))  - \lambda(X_{\varepsilon s}^N,\tilde Z_{\varepsilon s}^N,\tilde U_{\varepsilon s}^N(\cdot))\Big) ds \notag\\
	&-\int_0^{1} (Z_{t}^\varepsilon - \bar Z_{t}) d W^1_{t} -  \sqrt \varepsilon\int_0^{1}\Big (\check \Xi^{N}_{\frac{t}{\varepsilon}}-\hat \Xi^{\varepsilon}_\frac{t}{\varepsilon}\Big ) d   W^2_t\notag\\
	&-\int_0^{1}\int_{H \setminus 0} (U^{\varepsilon}_{t}(w)-\bar U_{t}(w)) (N_1(dt \, dw) - \nu_1(dw) dt)\notag\\
	&-  \varepsilon\int_0^{1}\int_{K \setminus 0} (\check\Theta^{N}_\frac{t}{\varepsilon}(w)-\hat\Theta^{\varepsilon}_\frac{t}{\varepsilon}(w)) \Big (N^\varepsilon_2(dt \, dw) - \frac{\nu_2(dw)}{\varepsilon} dt\Big ).
\end{align}
Now we notice that   
\begin{align*}
&\psi (X_{\varepsilon s}, \hat Q_s^{\varepsilon},Z^{\varepsilon}_{\varepsilon s},  \hat \Xi^{\varepsilon}_{s}, U^{\varepsilon}_{\varepsilon s}(\cdot), \hat \Theta^{\varepsilon}_s(\cdot) )-\psi (X_{\varepsilon s}, \hat Q_s^{\varepsilon},\bar Z_{\varepsilon s},  \hat \Xi^{\varepsilon}_s, \bar U_{\varepsilon s}(\cdot), \hat \Theta^{\varepsilon}_s(\cdot))\notag	\\
&+\psi(X_{\varepsilon s}, \hat Q_s^{\varepsilon},\bar Z_{\varepsilon s},  \hat \Xi^{\varepsilon}_{s}, \bar U_{\varepsilon s}(\cdot), \hat \Theta^{\varepsilon}_s(\cdot))-\psi(X^N_{\varepsilon s}, \hat {\mathcal{Q}}_s^{N}, \tilde Z_{\varepsilon s}^N, \check \Xi^{N}_s, \tilde U_{\varepsilon s}^N(\cdot), \check \Theta^{N}_s(\cdot))\\
%&=\psi\Big (X_{\varepsilon s}, \hat Q_s^{\varepsilon},Z^{\varepsilon}_{\varepsilon s},  \hat \Xi^{\varepsilon}_{s}, U^{\varepsilon}_{\varepsilon s}(\cdot), \hat \Theta^{\varepsilon}_s(\cdot)\Big )-\psi(X^N_{\varepsilon s}, \hat {\mathcal{Q}}_s^{N}, \tilde Z_{\varepsilon s}^N, \check \Xi^{N}_s, \tilde U_{\varepsilon s}^N(\cdot), \check \Theta^{N}_s(\cdot))\\
&= \psi(X_{\varepsilon s}, \hat Q_s^{\varepsilon},{Z^{\varepsilon}_{\varepsilon s}},  \hat \Xi^{\varepsilon}_{s}, U^{\varepsilon}_{\varepsilon s}(\cdot), \hat \Theta^{\varepsilon}_s(\cdot))-\psi (X_{\varepsilon s}, \hat Q_s^{\varepsilon},{\bar Z_{\varepsilon s}},  \hat \Xi^{\varepsilon}_s,  U_{\varepsilon s}^\varepsilon(\cdot), \hat \Theta^{\varepsilon}_s(\cdot))\notag	\\
&+ \psi(X_{\varepsilon s}, \hat Q_s^{\varepsilon},\bar Z_{\varepsilon s},  \hat \Xi^{\varepsilon}_{s}, {U^{\varepsilon}_{\varepsilon s}(\cdot)}, \hat \Theta^{\varepsilon}_s(\cdot))-\psi(X_{\varepsilon s}, \hat {Q}_s^{\varepsilon},\bar Z_{\varepsilon s},  \hat \Xi^{\varepsilon}_{s}, {\bar U_{\varepsilon s}(\cdot)}, \hat \Theta^{\varepsilon}_s(\cdot))\\
%&+\psi\Big (X_{\varepsilon s}, \hat { Q}_s^{\varepsilon},\bar Z_{\varepsilon s},  \hat \Xi^{\varepsilon}_{s}, \textcolor{red}{\bar U_{\varepsilon s}(\cdot)}, \hat \Theta^{\varepsilon}_s(\cdot)\Big )
%%&-\psi\Big (X_{\varepsilon s}^N, \hat {\mathcal Q}_s^{N},\bar Z_{\varepsilon s},  \hat \Xi^{\varepsilon}_{s}, \textcolor{red}{\tilde U^{N}_{\varepsilon s}(\cdot)}, \hat \Theta^{\varepsilon}_s(\cdot)\Big )\\
%%&+\psi\Big (X_{\varepsilon s}^N, \hat {\mathcal Q}_s^{N},\textcolor{red}{\bar Z_{\varepsilon s}},  \hat \Xi^{\varepsilon}_{s}, \tilde U^{N}_{\varepsilon s}(\cdot), \hat \Theta^{\varepsilon}_s(\cdot)\Big )\\
%-\psi\Big (X_{\varepsilon s}^N, \hat {\mathcal Q}_s^{N},\textcolor{red}{\tilde Z^N_{\varepsilon s}},  \hat \Xi^{\varepsilon}_{s}, \tilde U^{N}_{\varepsilon s}(\cdot), \hat \Theta^{\varepsilon}_s(\cdot)\Big )\\
&+R^{\varepsilon, N}_s\\
&+\psi (X_{\varepsilon s}^N, \hat {\mathcal Q}_s^{N},\tilde Z^N_{\varepsilon s},  {\hat \Xi^{\varepsilon}_{s}}, \tilde U^{N}_{\varepsilon s}(\cdot), \hat \Theta^{\varepsilon}_s(\cdot) )-\psi(X^N_{\varepsilon s}, \hat {\mathcal{Q}}_s^{N}, \tilde Z_{\varepsilon s}^N, {\check \Xi^{N}_s}, \tilde U_{\varepsilon s}^N(\cdot), \hat \Theta^{\varepsilon}_s(\cdot))\\
%-\psi\Big (X_{\varepsilon s}, \hat Q_s^{\varepsilon},\textcolor{blue}{\bar Z_{\varepsilon s}},  \hat \Xi^{\varepsilon}_s, \textcolor{blue}{\bar U_{\varepsilon s}(\cdot)}, \hat \Theta^{\varepsilon}_s(\cdot)\Big)\\
%+\psi\Big (X_{\varepsilon s}, \hat Q_s^{\varepsilon},\bar Z_{\varepsilon s},  \hat \Xi^{\varepsilon}_{s}, \bar U_{\varepsilon s}(\cdot), \hat \Theta^{\varepsilon}_s(\cdot)\Big )
&+\psi(X^N_{\varepsilon s}, \hat {\mathcal{Q}}_s^{N}, \tilde Z_{\varepsilon s}^N, \check \Xi^{N}_s, \tilde U_{\varepsilon s}^N(\cdot), {\hat \Theta^{\varepsilon}_s(\cdot)})-\psi(X^N_{\varepsilon s}, \hat {\mathcal{Q}}_s^{N}, \tilde Z_{\varepsilon s}^N, \check \Xi^{N}_s, \tilde U_{\varepsilon s}^N(\cdot), {\check \Theta^{N}_s(\cdot)}), 
\end{align*}
where we have set 
\begin{equation}\label{restoR}
R^{\varepsilon, N}_s := \psi({X_{\varepsilon s}}, {\hat {Q}_s^{\varepsilon}},{\bar Z_{\varepsilon s}},  \hat \Xi^{\varepsilon}_{s}, {\bar U_{\varepsilon s}(\cdot)}, \hat \Theta^{\varepsilon}_s(\cdot))
-\psi ({X_{\varepsilon s}^N}, {\hat {\mathcal Q}_s^{N}},{\tilde Z^N_{\varepsilon s}},  \hat \Xi^{\varepsilon}_{s}, {\tilde U^{N}_{\varepsilon s}(\cdot)}, \hat \Theta^{\varepsilon}_s(\cdot)).
\end{equation}
Recalling Hypothesis  
 {\sc{\textbf{(H$\psi$)}}}-c), there  exist $L_x, L_q, L_z, %L_\zeta, 
  L_u%, L_\theta
 >0$ such that 
	 \begin{align}\label{Reps_est}
		|R^{\varepsilon, N}_s|&\leq L_z|\bar Z_{\varepsilon s}-\tilde Z^N_{\varepsilon s}|_{\Xi^\ast}
		%+  L_\zeta|\zeta-\zeta'|_{\Xi^\ast}
		+ L_u \Big(\int_{H\setminus 0}|\bar U_{\varepsilon s}(w)- \tilde U^{N}_{\varepsilon s}(w)|^2 \nu_1(dw)\Big)^\frac{1}{2} 
		%+ L_\theta \Big(\int_{\Xi \setminus 0}|\theta_s(w)- \theta'_s(w)|^2 \nu_2(dw)\Big)^\frac{1}{2}
		\\
		&+ L_x \Big(1 + |\bar Z_{\varepsilon s}|_{\Xi^\ast} + \Big(\int_{H \setminus 0}|\bar U_{\varepsilon s}(w)|^2 \nu_1(dw)\Big)^\frac{1}{2}\Big)|X_{\varepsilon s}-X_{\varepsilon s}^N|_H\\
		&+L_q \Big(1 + |\bar Z_{\varepsilon s}|_{\Xi^\ast} + \Big(\int_{H \setminus 0}|\bar U_{\varepsilon s}(w)|^2 \nu_1(dw)\Big)^\frac{1}{2}\Big)|\hat {Q}_s^{\varepsilon}-\hat {\mathcal Q}_s^{N}|_K.
	\end{align}
At this point we introduce 
\begin{align*}
	\delta^{1,\varepsilon}(s) =
	\left\{
\begin{array}{ll}
 \frac{|\psi(X_{\varepsilon s}, \hat Q_s^{\varepsilon},{Z^{\varepsilon}_{\varepsilon s}},  \hat \Xi^{\varepsilon}_{s}, U^{\varepsilon}_{\varepsilon s}(\cdot), \hat \Theta^{\varepsilon}_s(\cdot))-\psi (X_{\varepsilon s}, \hat Q_s^{\varepsilon},{\bar Z_{\varepsilon s}},  \hat \Xi^{\varepsilon}_s,  U_{\varepsilon s}^\varepsilon(\cdot), \hat \Theta^{\varepsilon}_s(\cdot))|}{|Z^{\varepsilon}_{\varepsilon s}-\bar Z_{\varepsilon s}|^2}(Z^{\varepsilon}_{\varepsilon s}-\bar Z_{\varepsilon s})^\ast, & \textup{if}\,\,|Z^{\varepsilon}_{\varepsilon s}-\bar Z_{\varepsilon s}| \neq 0,\\
0, &\textup{otherwise}, 
\end{array}
\right. 
\end{align*}
\begin{align*}
	\delta^{2,\varepsilon, N}(s) =
	\left\{
\begin{array}{ll}
 \frac{|\psi (X_{\varepsilon s}^N, \hat {\mathcal Q}_s^{N},\tilde Z^N_{\varepsilon s},  {\hat \Xi^{\varepsilon}_{s}}, \tilde U^{N}_{\varepsilon s}(\cdot), \hat \Theta^{\varepsilon}_s(\cdot) )-\psi(X^N_{\varepsilon s}, \hat {\mathcal{Q}}_s^{N}, \tilde Z_{\varepsilon s}^N, {\check \Xi^{N}_s}, \tilde U_{\varepsilon s}^N(\cdot), \hat \Theta^{\varepsilon}_s(\cdot))|}{|\hat \Xi^{\varepsilon}_{s}-\check \Xi^{N}_s|^2}(\hat \Xi^{\varepsilon}_{s}-\check \Xi^{N}_s)^\ast, & \textup{if}\,\,|\hat \Xi^{\varepsilon}_{s}-\check \Xi^{N}_s| \neq 0,\\
0, &\textup{otherwise}.
\end{array}
\right. 
\end{align*}
We notice that, by Hypothesis {\sc{\textbf{(H$\psi$)}}}-c), we have the uniform bounds
\begin{equation}\label{unifbounddelta}
	|\delta^{1,\varepsilon}(s)| \leq L_z, \quad |\delta^{2,\varepsilon, N}(s)|\leq L_\zeta, \quad s \in \big[0, 1/\varepsilon\big ].
\end{equation}
%the processes $(\delta^{1,\varepsilon}(s))_{s \in [0, \frac{1}{\varepsilon}]}$ and $(\delta^{2,\varepsilon, N}(s))_{s \in [0, \frac{1}{\varepsilon}]}$ are bounded uniformly by $L_z$ and $L_\zeta$, respectively.
 Moreover,   
% $$
%f^{x,z,u}(q,\zeta, \theta) := \psi(x,q,z,\zeta,u, \theta), \quad x \in H, q \in K, z, \zeta \in \Xi^\ast, u \in L^2(\nu_1), \theta\in L^2(\nu_2).
%$$
%Then,
 by  Hypothesis {\sc{\textbf{(H$\psi$)}}}-d), Remark \ref{Cohencond} insures the existence of two  
 functions 
  $\gamma_s^{1, \varepsilon}$ and $\gamma_s^{2,\varepsilon,N}$, depending also on $U^\varepsilon, \bar{U},\hat{\Theta}^\varepsilon, \check{\Theta}^N$, such that
\begin{align}
	&C_1 (1 \wedge |w|_H) \leq \gamma^{1, \varepsilon}(w) \leq C_2 (1 \wedge |w|_\Xi), \quad C_1 \in (-1, 0], C_2 \geq 0,\label{est_nui}\\
	&\bar C_1 (1 \wedge |w|_K) \leq \gamma^{2,\varepsilon,N}(w) \leq \bar C_2 (1 \wedge |w|_K), \quad \bar C_1 \in (-1, 0], \bar C_2 \geq 0,\label{est_nui_bis}
\end{align}
   %:=\gamma^{\hat {\mathcal{Q}}_s^{N},  \check \Xi^{N}_s, \hat \Theta^{\varepsilon}_s(\cdot), \check \Theta^{N}_s(\cdot)}$ 
     and such that 
  		\begin{align}\label{est_f_gamma_3}
		&\psi(X_{\varepsilon s},\hat Q_s^{\varepsilon},\bar Z_{\varepsilon s},\hat \Xi^{\varepsilon}_{s},U^{\varepsilon}_{\varepsilon s}(\cdot),\hat \Theta^{\varepsilon}_s(\cdot))-\psi(X_{\varepsilon s},\hat Q_s^{\varepsilon},\bar Z_{\varepsilon s},\hat \Xi^{\varepsilon}_{s},\bar U_{\varepsilon s}(\cdot), \hat \Theta^{\varepsilon}_s(\cdot)) \\
		&= \int_{H \setminus 0} \gamma_s^{1, \varepsilon}(w)  (U^{\varepsilon}_{\varepsilon s}(w) -\bar U_{\varepsilon s}(w)) \nu_1(dw)
	\end{align}
	and 
\begin{align}\label{est_f_gamma_2}
		&\psi(X^N_{\varepsilon s},\hat {\mathcal{Q}}_s^{N}, \tilde Z_{\varepsilon s}^N,\check \Xi^{N}_s,\tilde U_{\varepsilon s}^N(\cdot), \hat \Theta^{\varepsilon}_s(\cdot))-\psi(X^N_{\varepsilon s},\hat {\mathcal{Q}}_s^{N}, \tilde Z_{\varepsilon s}^N,\check \Xi^{N}_s, \tilde U_{\varepsilon s}^N(\cdot),\check \Theta^{N}_s(\cdot))\\
		& = \int_{K \setminus 0} \gamma_s^{2, \varepsilon, N}(w)  (\hat \Theta^{\varepsilon}_s(w) -\check \Theta^{N}_s(w)) \nu_2(dw).
	\end{align}
		%\textcolor{red}{We have to verify that  
		% $\gamma^i: \Omega \times [0,\,T] \times \Xi\rightarrow \R_+$ are  uniformly bounded functions such that   $\gamma^i_t(\omega, \cdot) -1\in L^2(\nu^i)$ for all $(\omega, t) \in \Omega \times [0,\,T]$
		% the terms in $\nu^i$ are well defined (OK if $\gamma^i$ are uniformly bounded).}		
Identity \eqref{diff_to_est_2} reads 
\begin{align}\label{diff_to_est_3}
	Y^\varepsilon_0 -\bar Y_0&=\varepsilon\sum_{k=0}^{2^N-1}\Big  (\check Y_{\frac{t_k}{\varepsilon}}^{N,k} - \check Y^{N,k}_\frac{t_{k+1}}{\varepsilon}\Big )+\varepsilon\int_0^{\frac{1}{\varepsilon}}R^{\varepsilon, N}_s  ds\notag\\
	&+\varepsilon\int_0^{\frac{1}{\varepsilon}} 
	\delta^{1, \varepsilon}_s (Z^{\varepsilon}_{\varepsilon s}-\bar Z_{\varepsilon s}) ds+\varepsilon\int_0^{\frac{1}{\varepsilon}}\delta^{2, \varepsilon, N}_s(\hat \Xi^{\varepsilon}_{s}-\check \Xi^{N}_s) ds\notag\\
	&+\varepsilon\int_0^{\frac{1}{\varepsilon}}\int_{H \setminus 0} \gamma_s^{1, \varepsilon}(w)  (U^{\varepsilon}_{\varepsilon s}(w) -\bar U_{\varepsilon s}(w)) \nu_1(dw) ds\notag\\
	&+\varepsilon\int_0^{\frac{1}{\varepsilon}}\int_{K \setminus 0} \gamma_s^{2, \varepsilon, N}(w)  (\hat \Theta^{\varepsilon}_s(w) -\check \Theta^{N}_s(w)) \nu_2(dw) ds\notag\\
	& - \varepsilon\int_0^{\frac{1}{\varepsilon}} \Big(\lambda(X_{\varepsilon s},\bar Z_{\varepsilon s},\bar U_{\varepsilon s}(\cdot))  - \lambda(X_{\varepsilon s}^N,\tilde Z_{\varepsilon s}^N,\tilde U_{\varepsilon s}^N(\cdot))\Big) ds \notag\\
	&-\int_0^{1} (Z_{t}^\varepsilon - \bar Z_{t}) d W^1_{t} -  \sqrt \varepsilon\int_0^{1}\Big (\check \Xi^{N}_{\frac{t}{\varepsilon}}-\hat \Xi^{\varepsilon}_\frac{t}{\varepsilon}\Big ) d   W^2_t\notag\\
	&-\int_0^{1}\int_{H \setminus 0} (U^{\varepsilon}_{t}(w)-\bar U_{t}(w)) (N_1(dt \, dw) - \nu_1(dw) dt)\notag\\
	&-  \varepsilon\int_0^{1}\int_{K \setminus 0} (\check\Theta^{N}_\frac{t}{\varepsilon}(w)-\hat\Theta^{\varepsilon}_\frac{t}{\varepsilon}(w)) \Big (N^\varepsilon_2(dt \, dw) - \frac{\nu_2(dw)}{\varepsilon} dt\Big ), 
\end{align}
and, rescaling the time, it gives
	 \begin{align}\label{diff_to_est_4}
	Y^\varepsilon_0 -\bar Y_0&=\varepsilon\sum_{k=0}^{2^N-1}\Big  (\check Y_{\frac{t_k}{\varepsilon}}^{N,k} - \check Y^{N,k}_\frac{t_{k+1}}{\varepsilon}\Big )+\int_0^{1}R^{\varepsilon, N}_{\frac{t}{\varepsilon}}  dt\notag\\
	& - \int_0^{1} \Big(\lambda(X_{t},\bar Z_{t},\bar U_{t}(\cdot))  - \lambda(X_{t}^N,\tilde Z_{t}^N,\tilde U_{t}^N(\cdot))\Big) dt \notag\\
	&+\int_0^{1} 
	\delta^{1, \varepsilon}_{\frac{t}{\varepsilon}} (Z^{\varepsilon}_{t}-\bar Z_{t}) dt-\int_0^{1} (Z_{t}^\varepsilon - \bar Z_{t}) d W^1_{t} \notag\\
	&+\int_0^{1}\delta^{2, \varepsilon, N}_{\frac{t}{\varepsilon}}\Big (\hat \Xi^{\varepsilon}_{\frac{t}{\varepsilon}}-\check \Xi^{N}_{\frac{t}{\varepsilon}}\Big ) dt-  \sqrt \varepsilon\int_0^{1}\Big (\check \Xi^{N}_{\frac{t}{\varepsilon}}-\hat \Xi^{\varepsilon}_\frac{t}{\varepsilon}\Big ) d   W^2_t\notag\\
	&-\int_0^{1}\int_{H \setminus 0} (U^{\varepsilon}_{t}(w)-\bar U_{t}(w))\Big  (N_1(dt \, dw) -  {\big (\gamma_{\frac{t}{\varepsilon}}^{1, \varepsilon}(w)+1\big )}\nu_1(dw) dt\Big )\\
	%&+\int_0^{1}\int_{\Xi \setminus 0}  (U^{\varepsilon}_{t}(w) -\bar U_{t}(w)) \nu_1(dw) dt\notag\\
	&-  \varepsilon\int_0^{1}\int_{K \setminus 0} \Big (\check\Theta^{N}_\frac{t}{\varepsilon}(w)-\hat\Theta^{\varepsilon}_\frac{t}{\varepsilon}(w)\Big ) \Big (N^\varepsilon_2(dt \, dw) - {\big (\gamma_{\frac{t}{\varepsilon}}^{2, \varepsilon, N}(w)+1\big )}\frac{\nu_2(dw)}{\varepsilon} dt\Big ).\notag
	%&+\int_0^{1}\int_{\Xi \setminus 0}  \Big (\hat \Theta^{\varepsilon}_{\frac{t}{\varepsilon}}(w) -\check \Theta^{N}_{\frac{t}{\varepsilon}}(w)\Big ) \nu_2(dw) dt.\notag
\end{align}
We set, for $s \in [0,1]$,   
\begin{align}
	M^1_s&:=\int_{[0,\,s]}\delta^{1, \varepsilon}_{\frac{t}{\varepsilon}} d W^1_t + \int_{]0,\,s]} \int_{H \setminus 0}{\gamma^{1, \varepsilon}_{\frac{t}{\varepsilon}}(w)} (N_1(dt\,dw)-\nu_1(dw)dt),\label{M1}\\ 
	M^2_s&:=\int_{[0,\,s]}\frac{1}{\sqrt \varepsilon  }\delta^{2, \varepsilon, N}_{\frac{t}{\varepsilon}} d W^2_t + \int_{]0,\,s]} \int_{K \setminus 0}{\gamma^{2, \varepsilon, N}_{\frac{t}{\varepsilon}}(w)} \Big (N^{\varepsilon}_2(dt\,dw)-\frac{\nu_2(dw)}{\varepsilon}dt\Big ).\label{M2}
\end{align}
Notice that  the jumps of  $M^1$ and $M^2$ in \eqref{M1}-\eqref{M2} are completely described by the discontinuities of the processes $X$ and $Q$ (see system \eqref{fbsystem_2}), namely \begin{equation*}
	\Delta M^1_s =  \gamma^{1, \varepsilon}_{\frac{s}{\varepsilon}}(\Delta X_{s}), \qquad \Delta M^2_s:=  \gamma^{2, \varepsilon, N}_{\frac{s}{\varepsilon}}(\Delta Q^\varepsilon_{\frac{s}{\varepsilon}}), \quad s\in [0,1],
\end{equation*}
where  for any càdlàg process $L$ we have set
\begin{equation}\label{salto def}
	\Delta L_s:=  L_s - L_{s^-}.  
\end{equation}
We define (see \eqref{DDexp}) 
\begin{align*}
	\frac{d \P^{\varepsilon, N}}{d\P}:= \mathcal E_1(M^1+M^2)
	&=e^{\int_0^1\delta^{1, \varepsilon}_{\frac{t}{\varepsilon}} d W^1_t-\frac{1}{2} \int_0^1\big |\delta^{1, \varepsilon}_{\frac{t}{\varepsilon}}\big |^2 dt}e^{-\int_{]0,\,1]} \int_{H \setminus 0} \gamma_{\frac{t}{\varepsilon}}^{1, \varepsilon}(\ell) \nu_1(d\ell)dt}\\
	&\cdot e^{\int_0^1\frac{1}{\sqrt \varepsilon} \delta^{2, \varepsilon, N}_{\frac{t}{\varepsilon}} d W^2_t-\frac{1}{2} \int_0^1\frac{1}{ \varepsilon}\big |\delta^{2, \varepsilon, N}_{\frac{t}{\varepsilon}}\big |^2 dt}e^{-\int_{]0,\,1]} \int_{K \setminus 0}\gamma_{\frac{t}{\varepsilon}}^{2, \varepsilon, N}(\ell) \frac{1}{\varepsilon}\nu_2(d\ell)dt}\\
	&\cdot \prod_{n \geq 1: T_n, S_n \leq 1}\Big(\gamma^{1, \varepsilon}_{\frac{T_n}{\varepsilon}}(\Delta X_{T_n})+\gamma^{2, \varepsilon, N}_{\frac{S_n}{\varepsilon}}(\Delta Q^\varepsilon_{\frac{S_n}{\varepsilon}})+1\Big),
\end{align*}
where  $T_n$ (resp. $S_n$) denote the jump times of $N_1(dt\,dw)$ (resp. $N_2(dt\,dw)$).

By the Girsanov Theorem \ref{Girsanov} {(notice that by \eqref{est_nui}  $\gamma^{i,\varepsilon}$ is uniformly bounded and   $\gamma^{i, \varepsilon}(\omega, t, \cdot) \in L^2(\nu_i)$, $i =1,2$,  for all $(\omega, t) \in \Omega \times [0,\,T]$)}, the compensator of $(N_1(dt\,dz),N^\varepsilon_2(dt\,dz))$ under $\P^{\varepsilon}$ is 
$$
{\Big(\big (\gamma^{1, \varepsilon}_{\frac{t}{\varepsilon}}(w)+1\big )}\nu_1(dw) dt, {\big (\gamma^{2, \varepsilon, N}_{\frac{t}{\varepsilon}}(w)+1\big )}\frac{\nu_2(dw)}{\varepsilon} dt\Big), 
$$
 and the process
	\begin{align*}
	(\tilde W^1,\tilde W^2)&:=\Big (W^1 - \int_{[0,\,\cdot]}\delta^{1, \varepsilon}_{\frac{t}{\varepsilon}} dt, W^2 - \int_{[0,\,\cdot]}\frac{1}{\sqrt \varepsilon}   \delta^{2, \varepsilon, N}_{\frac{t}{\varepsilon}}dt\Big)
	%(N_1^{\P^{\varepsilon, \alpha}}(dt\,dz),N_2^{\varepsilon,\P^{\varepsilon, \alpha}}(dt\,dz))&:=,
\end{align*}
is such that {$\tilde W^1$ and $\tilde W^2$ are cylindrical Wiener  processes under $\P^{\varepsilon, N}$.} 
 Taking 	 the expectation ${\E}^{\varepsilon, N}$ under the new probability $\P^{\varepsilon, N}$ in \eqref{diff_to_est_4}, we get 
	 \begin{align}\label{diff_to_est_5}
	Y^\varepsilon_0 -\bar Y_0&=\varepsilon\sum_{k=0}^{2^N-1}{\E}^{\varepsilon, N}\Big  [\check Y_{\frac{t_k}{\varepsilon}}^{N,k} - \check Y^{N,k}_\frac{t_{k+1}}{\varepsilon}\Big ]+{\E}^{\varepsilon, N} \Big [\int_0^{1}R^{\varepsilon, N}_{\frac{t}{\varepsilon}}  dt\Big ]\notag\\
%	&+\int_0^{1} 
%	\delta^{1, \varepsilon}_{\frac{t}{\varepsilon}} (Z^{\varepsilon}_{t}-\bar Z_{t}) dt-\int_0^{1} (Z_{t}^\varepsilon - \bar Z_{t}) d W^1_{t} \notag\\
	& -{\E}^{\varepsilon, N} \Big[ \int_0^{1} \Big(\lambda(X_{t},\bar Z_{t},\bar U_{t}(\cdot))  - \lambda(X_{t}^N,\tilde Z_{t}^N,\tilde U_{t}^N(\cdot))\Big) dt \Big ].\notag
	%&+\int_0^{1}\delta^{2, \varepsilon, N}_{\frac{t}{\varepsilon}}\Big (\hat \Xi^{\varepsilon}_{\frac{t}{\varepsilon}}-\check \Xi^{N}_{\frac{t}{\varepsilon}}\Big ) dt-  \sqrt \varepsilon\int_0^{1}\Big (\check \Xi^{N}_{\frac{t}{\varepsilon}}-\hat \Xi^{\varepsilon}_\frac{t}{\varepsilon}\Big ) d   W^2_t\notag\\
%	&-\int_0^{1}\int_{\Xi \setminus 0} (U^{\varepsilon}_{t}(w)-\bar U_{t}(w))\Big  (N_1(dt \, dw) -  \gamma_{\frac{t}{\varepsilon}}^{1, \varepsilon}(w)\nu_1(dw) dt\Big )\notag\\
%	&+{\E}^{\varepsilon, N} \Big[\int_0^{1}\int_{\Xi \setminus 0}  (U^{\varepsilon}_{t}(w) -\bar U_{t}(w)) \nu_1(dw) dt\Big ]\notag\\
%	&-  \varepsilon\int_0^{1}\int_{\Xi \setminus 0} \Big (\check\Theta^{N}_\frac{t}{\varepsilon}(w)-\hat\Theta^{\varepsilon}_\frac{t}{\varepsilon}(w)\Big ) \Big (N^\varepsilon_2(dt \, dw) - \gamma_{\frac{t}{\varepsilon}}^{2, \varepsilon, N}(w)\frac{\nu_2(dw)}{\varepsilon} dt\Big )\\
%	&+{\E}^{\varepsilon, N} \Big[\int_0^{1}\int_{\Xi \setminus 0}  \Big (\hat \Theta^{\varepsilon}_{\frac{t}{\varepsilon}}(w) -\check \Theta^{N}_{\frac{t}{\varepsilon}}(w)\Big ) \nu_2(dw) dt\Big ].\notag
	\end{align}
	Recalling estimate  \eqref{Reps_est}, and the fact that, by Theorem \ref{T:ergodic},
there exist some positive constants $L'_x$, $L'_z$, $L'_u$ such that 
\begin{align*}
	&|\lambda(X_{t},\bar Z_{t},\bar U_{t}(\cdot))  - \lambda(X_{t}^N,\tilde Z_{t}^N,\tilde U_{t}^N(\cdot))|\\
	&\leq L'_x\Big(1+ |\bar Z_{t}|_{\Xi^\ast} + \Big(\int_{H\setminus 0}|\bar U_{t}(w)|^2 \nu_1(dw)\Big)^\frac{1}{2}\Big)|X_{t}-X_{t}^N|_H\notag\\
	&+L'_z |\bar Z_{t}-\tilde Z_{t}^N|_{\Xi^\ast}+ L'_u\Big(\int_{H \setminus 0}|\bar U_{t}(w)- \tilde U_{t}^N(w)|^2 \nu_1(dw)\Big)^\frac{1}{2}, 
\end{align*}
   previous identity yields 
	\begin{align}
	&Y^\varepsilon_0 -\bar Y_0\leq \varepsilon\sum_{k=0}^{2^N-1}{\E}^{\varepsilon, N}\Big  [\check Y_{\frac{t_k}{\varepsilon}}^{N,k} - \check Y^{N,k}_\frac{t_{k+1}}{\varepsilon}\Big ]+(L_z+ L'_z){\E}^{\varepsilon, N} \Big [\int_0^{1}|\bar Z_{t}-\tilde Z^N_{t}|_{\Xi^\ast}dt\Big ]
		%+  L_\zeta|\zeta-\zeta'|_{\Xi^\ast}
		\notag\\
		&+ (L_u+L'_u) {\E}^{\varepsilon, N} \Big [\int_0^{1}\Big(\int_{H \setminus 0}|\bar U_{t}(w)- \tilde U^{N}_{t}(w)|^2 \nu_1(dw)\Big)^\frac{1}{2}dt \Big]
		%+ L_\theta \Big(\int_{\Xi \setminus 0}|\theta_s(w)- \theta'_s(w)|^2 \nu_2(dw)\Big)^\frac{1}{2}
		\notag\\
		&+ (L_x+L'_x) {\E}^{\varepsilon, N} \Big [\int_0^{1}\Big (1 + |\bar Z_{t}|_{\Xi^\ast}+\Big(\int_{H \setminus 0}|\bar U_{t}(w)|^2 \nu_1(dw)\Big)^\frac{1}{2}\Big )|X_{t}-X_{t}^N|_H dt\Big]\\
		%&+ (L_x+L'_x) {\E}^{\varepsilon, N} \Big [\int_0^{1}\Big(\int_{\Xi \setminus 0}|\bar U_{t}(w)|^2 \nu_1(dw)\Big)^\frac{1}{2}|X_{t}-X_{t}^N|_H dt\Big]\notag\\
		&
		+L_q {\E}^{\varepsilon, N} \Big [\int_0^{1}\Big (1 + |\bar Z_{t}|_{\Xi^\ast}+\Big(\int_{H \setminus 0}|\bar U_{t}(w)|^2 \nu_1(dw)\Big)^\frac{1}{2}\Big )\big|\hat {Q}_{\frac{t}{\varepsilon}}^{\varepsilon}-\hat {\mathcal Q}_{\frac{t}{\varepsilon}}^{N}\big |_K  dt\Big ]
		\notag\\
			%	&
		%+L_q {\E}^{\varepsilon, N} \Big [\int_0^{1} \Big(\int_{\Xi \setminus 0}|\bar U_{t}(w)|^2 \nu_1(dw)\Big)^\frac{1}{2}\big|\hat {Q}_{\frac{t}{\varepsilon}}^{\varepsilon}-\hat {\mathcal Q}_{\frac{t}{\varepsilon}}^{N}\big |_K  dt\Big ]\notag\\
		&=: I_1 + I_2+ I_3 + I_4 + I_5.\label{lastid}
\end{align}
We notice that 
\begin{align}\label{F-B_limit_system_3}
\left\{
\begin{array}{llll}
&d X_t=   A X_t dt - R \delta^{1, \varepsilon}_{\frac{t}{\varepsilon}}dt
%+  R\int_{\Xi \setminus 0} w \gamma_{\frac{t}{\varepsilon}}^{1, \varepsilon}(w)\nu_1(dw) dt
 + R d \tilde W_t^1 + R\int_{H \setminus 0} w [N_1(dt\,dw) - %( \gamma_{\frac{t}{\varepsilon}}^{1, \varepsilon}(w)+1)
\nu_1(dw) dt], \\
&X_0 = x_0,\\
& -d\bar Y_t =    \lambda(X_t,\bar Z_t,\bar U_t(\cdot)) dt  -  \bar Z_t \big [d  \tilde W^1_t- \delta^{1, \varepsilon}_{\frac{t}{\varepsilon}}dt\big ]
% -  \int_{\Xi \setminus 0} \bar U_t(w)  \gamma_{\frac{t}{\varepsilon}}^{1, \varepsilon}(w)\nu_1(dw) dt
-  \int_{H \setminus 0} \bar U_t(w) [N_1(dt \, dw) - %( \gamma_{\frac{t}{\varepsilon}}^{1, \varepsilon}(w)+1)
	  \nu_1(dw) dt], \\
&	\bar Y_1=h(X_1). 
\end{array}
\right.
\end{align}
We apply the Girsanov Theorem \ref{Girsanov_gen}  with $\nu(\omega, dt\,d\ell) =\big (\gamma^{1, \varepsilon}_{\frac{t}{\varepsilon}}(\omega, \ell)+1\big )\nu_1(d\ell) dt$, $\beta(\omega, t)= \delta^{1,\varepsilon}_\frac{t}{\varepsilon}(\omega)$ and  $\Gamma(\omega, t, \ell)= \big (\gamma^{1, \varepsilon}_{\frac{t}{\varepsilon}}(\omega, \ell)+1\big )^{-1}$.   
	If follows that,  if $ \tilde{\P}^\varepsilon:= \rho \, \P^{\varepsilon, N}$,  with 
\begin{align}\label{rho}
		 \rho 
&=e^{\int_0^1\delta^{1,\varepsilon}_\frac{t}{\varepsilon}  d \tilde W_t^1-\frac{1}{2} \int_0^1|\delta^{1,\varepsilon}_\frac{t}{\varepsilon}|^2 dt}e^{\int_{]0,\,1]} \int_{H \setminus 0}\gamma^{1, \varepsilon}_{\frac{t}{\varepsilon}}(w)\nu_1(dw)dt}\prod_{n \geq 1: T_n \leq 1}\big (\gamma^{1, \varepsilon}_{\frac{T_n}{\varepsilon}}(\Delta X_{T_n})+1\big )^{-1}, 
	\end{align}
%$$
%	M:=\int_{[0,\,\cdot]}\delta^{1,\varepsilon}_\frac{t}{\varepsilon} d W_t - \int_{]0,\,\cdot]} \int_{\Xi \setminus 0}\frac{\gamma^{1, \varepsilon}_{\frac{t}{\varepsilon}}(\omega, w)}{\gamma^{1, \varepsilon}_{\frac{t}{\varepsilon}}(\omega, w)+1} (N(dt\,dw)-\textcolor{blue}{\big (\gamma^{1, \varepsilon}_{\frac{t}{\varepsilon}}(\omega, w)+1\big )}\nu(dw)dt). 
%	$$
	then 
	$$
	 \tilde W^1_s - \int_{[0,\,s]}\delta^{1,\varepsilon}_\frac{t}{\varepsilon} dt, \quad s \in [0,1], 
	$$
	is a  cylindrical Wiener process, and the compensator of $N_1(\omega, dt\,d\ell)$   is $\nu_1(d\ell) dt$, {while $\tilde{W}_2$ remains a Wiener process and the compensator of $ N_2^\varepsilon(\omega, dt\,d\ell)$   is $\displaystyle(\gamma^{2,\varepsilon,N}_\frac{t}{\varepsilon} (\omega, \ell) +1) \frac{\nu_2(d\ell)}{\varepsilon}\, dt$.

In particular, by the uniqueness of the solution of the forward backward system \eqref{F-B_limit_system_3} (see Theorem \ref{T_wellposFB}), 
the law of the process $(X_t)_{t \geq 0}$ under $\tilde{\P^{\varepsilon}}$ coincides with its law under $\P$. Moreover, also the laws of $(\bar Z_t)_{t \geq 0}$ and $(\tilde Z^N_t)_{t \geq 0}$ and $(\bar U_t(\cdot))_{t \geq 0}$ and $(\tilde U^N_t(\cdot))_{t \geq 0}$ under $\tilde{\P}^{\varepsilon}$ coincide with the corresponding laws under $\P$. As a matter of fact, again by  Theorem \ref{T_wellposFB},  $\bar Z_t= \bar \zeta(X_t)$ (resp.  
	$\bar U_t(\cdot)= \bar \theta(X_t)(\cdot)$) where  $\bar \zeta$ (resp. $\bar \theta$) is a deterministic Borel-measurable function from $H$ to $\Xi^\ast$ (resp. to $L^2(\nu_1)$), so the laws of $(\bar Z_t)_{t \geq 0}$ and $(\tilde Z^N_t)_{t \geq 0}$ (resp. of $(\bar U_t(\cdot))_{t \geq 0}$ and $(\tilde U^N_t(\cdot))_{t \geq 0}$) depend only on the law of $X$ in a non-anticipating way. From these considerations it follows that 
	\begin{align}
		\E^{\varepsilon, N}\Big[\rho \int_0^1 |\bar Z_t|_{\Xi^\ast}^2 dt\Big]&=\E\Big[ \int_0^1 |\bar Z_t|_{\Xi^\ast}^2 dt\Big]< + \infty,\label{rhoZ} \\
		\E^{\varepsilon, N}\Big[\rho \int_0^1 \int_{H\setminus 0} |\bar U_t(w)|^2 \nu_1(dw) dt\Big]&=\E\Big[ \int_0^1 \int_{H\setminus 0} |\bar U_t(w)|^2 \nu_1(dw) dt\Big]< + \infty.\label{rhoU}
	\end{align}
%CONSIDERAZIONI VECCHIE 

Now we notice that, by \eqref{unifbounddelta} and \eqref{est_nui}, there exist constants $\eta$, $C$ and  $M$ such that   
%\textcolor{red}{CONTROLLARE SE LE COSTANTI SONO QUELLE PIU' COMODE} 
\begin{equation} \label{stime unif girsanov}
|\gamma^{1,\varepsilon}_\frac{t}{\varepsilon}(w)| \leq  C_1 \vee C_2, \qquad  0 < \eta \leq (\gamma^{1,\varepsilon}_\frac{t}{\varepsilon}(w)+1)^{-1} \leq C, \quad |\delta^{1,\varepsilon}_{t/\varepsilon}| \leq M.
\end{equation}
In particular , for every $p \geq 1$, 
$$
(\gamma^{1,\varepsilon}_\frac{t}{\varepsilon}(w)+1)^{p}\leq \frac{1}{\eta^p}.
$$
Recalling \eqref{rho}, for every $p \geq 1$  we have
%\begin{align*}
%		 \rho 
%		 &=e^{\int_0^1\delta^{1,\varepsilon}_\frac{t}{\varepsilon}  d \tilde W_t-\frac{1}{2} \int_0^1|\delta^{1,\varepsilon}_\frac{t}{\varepsilon}|^2 dt}e^{\int_{]0,\,1]} \int_{\Xi \setminus 0}\gamma^{1, \varepsilon}_{\frac{t}{\varepsilon}}(w)\nu_1(dw)dt}\prod_{n \geq 1: T_n \leq 1}\big (\gamma^{1, \varepsilon}_{\frac{T_n}{\varepsilon}}(X_{T_N})+1\big )^{-1}. 
%	\end{align*}
%	and therefore 
	\begin{align*}
		\rho^{-p}&=e^{-\int_0^1p\,\delta^{1,\varepsilon}_\frac{t}{\varepsilon}  d \tilde W_t+\frac{1}{2} \int_0^1 p|\delta^{1,\varepsilon}_\frac{t}{\varepsilon}|^2 dt}e^{-\int_{]0,\,1]} \int_{H \setminus 0} p \,\gamma^{1, \varepsilon}_{\frac{t}{\varepsilon}}(w)\nu_1(dw)dt}{\prod_{n \geq 1: T_n \leq 1}\big (\gamma^{1, \varepsilon}_{\frac{T_n}{\varepsilon}}(\Delta X_{T_n})+1\big )^{p}}\\
		&=e^{-\int_0^1p\,\delta^{1,\varepsilon}_\frac{t}{\varepsilon}  d \tilde W_t+\frac{1}{2} \int_0^1 p|\delta^{1,\varepsilon}_\frac{t}{\varepsilon}|^2 dt}e^{-\int_{]0,\,1]} \int_{H \setminus 0} p \,\gamma^{1, \varepsilon}_{\frac{t}{\varepsilon}}(w)\nu_1(dw)dt}\\
		&e^{\pm \int_{]0,\,1]} \int_{H \setminus 0}\big [1-\big (\gamma^{1, \varepsilon}_{\frac{t}{\varepsilon}}(w)+1\big )^p\big ] \big (\gamma^{1, \varepsilon}_{\frac{t}{\varepsilon}}(w)+1\big )\nu_1(d w)dt} e^{\pm \frac{1}{2} \int_0^1 p^2|\delta^{1,\varepsilon}_\frac{t}{\varepsilon}|^2 dt}{\prod_{n \geq 1: T_n \leq 1}\big (\gamma^{1, \varepsilon}_{\frac{T_n}{\varepsilon}}(\Delta X_{T_n})+1\big )^{p}}\\
		&=e^{\frac{1}{2} \int_0^1 (p+p^2)|\delta^{1,\varepsilon}_\frac{t}{\varepsilon}|^2 dt}e^{\int_{]0,\,1]} \int_{H \setminus 0}\big (-1-\gamma^{1, \varepsilon}_{\frac{t}{\varepsilon}}(w)(p+1) +\big (\gamma^{1, \varepsilon}_{\frac{t}{\varepsilon}}(w)+1\big )^{p+1}\big )\nu_1(dw)dt}\tilde \rho_p(1) , 
	\end{align*}
with $\tilde \rho_p$  the exponential martingale defined as 
\begin{multline} 
\label{expmart}
\tilde {\mathcal\rho_p}(r):= e^{-\int_0^r p \, \delta^{1,\varepsilon}_\frac{t}{\varepsilon} d \tilde W_t-\frac{1}{2} \int_0^r p^2|\delta^{1,\varepsilon}_\frac{t}{\varepsilon}|^2 dt}e^{\int_{]0,\,r]} \int_{H \setminus 0}\big [1-\big (\gamma^{1, \varepsilon}_{\frac{t}{\varepsilon}}(w)+1\big )^p\big ] \big (\gamma^{1, \varepsilon}_{\frac{t}{\varepsilon}}(w)+1\big )\nu_1(dw)dt} \\ \cdot {\prod_{n \geq r: T_n \leq 1}\big (\gamma^{1, \varepsilon}_{\frac{T_n}{\varepsilon}}(\Delta X_{T_n})+1\big )^{p}}.
\end{multline}
It follows that for every $p \geq 1$ 
\begin{equation}\label{stima_path}
		\rho^{-p}\leq  \kappa_p\, \tilde {\rho}_p(1),
\end{equation}
%for a suitable $\kappa_p >0$.}
where  we have set
\begin{equation}
 \kappa_p:=e^{\frac{1}{2}  (p-p^2)M^2}e^{(p\big(1-\frac{1}{C}\big)-\frac{1}{C}+\frac{1}{\eta^{p+1}})\nu_1(H)}.
 \end{equation}
Since $\E^{\varepsilon,N}[\tilde {\rho}_p(1)]=1$, we get that, for every $p \geq 1$, 
\begin{equation}\label{est_kappa}
\E^{\varepsilon,N}[\rho^{-p}] \leq \kappa_p.
\end{equation}

Let us now go back to \eqref{lastid}. We have to estimate the terms $I_i$, $i =1,..., 7$.   We set 
\begin{align*}
	&\Delta_{X, N}:= \sup_{t \in [0,1]} |X_t- X_t^N|_H, \quad \Delta_{Z, N}:= \int_0^1 |\bar Z_t- \tilde Z_t^N|_\Xi^2  dt,\\
	& \Delta_{U, N}:=\int_0^1 \int_{H \setminus 0}|\bar U_{t}(w)- \tilde U^{N}_{t}(w)|^2 \nu_1(dw)dt, 
	% \quad \Delta_{Q, N, \varepsilon}:= \sup_{t \in [0,1]} |\hat Q^\varepsilon_t- \hat {\mathcal Q}_t^N|_K,
\end{align*}
	where $X^N$, $\tilde Z^N$ and  $\tilde U^N$ %and  $\hat {\mathcal Q}^N$ 
	are the processes defined in \eqref{XN}-\eqref{ZN}-\eqref{UN}. %- \eqref{mathcalQ^N}.

Below $C$ will denote a constant independent of $N$ and $\varepsilon$, that may vary from line to line. 
 Let us start by considering  $I_4$. Applying the H\"older inequality and using  \eqref{rhoZ}-\eqref{rhoU} and \eqref{est_kappa}, we have 
 \begin{align}\label{estI4}
	I_4 &= (L_x + L'_{x}){\E}^{\varepsilon, N} \Big [\int_0^{1}\Big (1 + |\bar Z_{t}|_{\Xi^\ast}+\Big(\int_{H \setminus 0}|\bar U_{t}(w)|^2 \nu_1(dw)\Big)^\frac{1}{2}\Big )|X_{t}-X_{t}^N|_H dt\Big]\notag\\
	&\leq (L_x + L'_{x}){\E}^{\varepsilon, N} \Big [\Delta_{X, N}\int_0^{1}\Big (1 + |\bar Z_{t}|_{\Xi^\ast}+\Big(\int_{H \setminus 0}|\bar U_{t}(w)|^2 \nu_1(dw)\Big)^\frac{1}{2}\Big ) dt\Big]\notag\\
	&= (L_x + L'_{x}){\E}^{\varepsilon, N} \Big [\rho^{-3/4}(\rho^{1/4}\Delta_{X, N})\rho^{1/2}\int_0^{1}\Big (1 + |\bar Z_{t}|_{\Xi^\ast}+\Big(\int_{H \setminus 0}|\bar U_{t}(w)|^2 \nu_1(dw)\Big)^\frac{1}{2}\Big ) dt\Big]\notag\\
	&\leq (L_x + L'_{x})({\E}^{\varepsilon, N} [\rho^{-3}])^{1/4}({\E}^{\varepsilon, N}[\rho\,\Delta_{X, N}^4])^{1/4} \cdot \notag\\
	&\cdot \Big(2{\E}^{\varepsilon, N}\Big [\rho \int_0^{1} \Big (1 + |\bar Z_{t}|^2_{\Xi^\ast}+\int_{H \setminus 0}|\bar U_{t}(w)|^2 \nu_1(dw)\Big ) dt\Big]\Big)^{1/4}\notag\\
	& \leq C \,({\E}[\Delta_{X, N}^4])^{1/4}. 
\end{align}
Concerning $I_2$, applying the H\"older inequality,  \eqref{rhoZ}  for $\bar Z_t - \tilde Z^N_t$ in place of $\bar Z_t$,  and \eqref{est_kappa}, we get 
\begin{align}\label{I2est}
	I_2 & = (L_z+ L'_z){\E}^{\varepsilon, N} \Big [\int_0^{1}|\bar Z_{t}-\tilde Z^N_{t}|_{\Xi^\ast}dt\Big ] \notag\\
	&= (L_z+ L'_z){\E}^{\varepsilon, N} \Big [\rho^{-1/2} \rho^{1/2} \int_0^{1}|\bar Z_{t}-\tilde Z^N_{t}|_{\Xi^\ast}dt\Big ] \notag\\
	&\leq (L_z+ L'_z)({\E}^{\varepsilon, N} [\rho^{-1}])^{1/2}\Big({\E}^{\varepsilon, N} \Big [\rho\int_0^{1}|\bar Z_{t}-\tilde Z^N_{t}|^2_{\Xi^\ast}dt\Big ]\Big)^{1/2}\notag\\
	&\leq  C({\E} [\Delta_{Z, N}])^{1/2}. 
\end{align}
Analogously, for $I_3$ (using \eqref{rhoU}  for $\bar U_t - \tilde U^N_t$ in place of $\bar U_t$) we have 
\begin{align}\label{estI3}
I_3&= (L_u+L'_u) {\E}^{\varepsilon, N} \Big [\int_0^{1}\Big(\int_{H \setminus 0}|\bar U_{t}(w)- \tilde U^{N}_{t}(w)|^2 \nu_1(dw)\Big)^\frac{1}{2}dt \Big]\notag\\
&\leq (L_z+ L'_z)({\E}^{\varepsilon, N} [\rho^{-1}])^{1/2}\Big({\E}^{\varepsilon, N} \Big [\rho\int_0^1 \int_{H\setminus 0}|\bar U_{t}(w)- \tilde U^{N}_{t}(w)|^2 \nu_1(dw)dt \Big]\Big)^{1/2}\notag\\
&\leq C ({\E} [\Delta_{U, N}])^{1/2}.
\end{align}
Moreover
\begin{align}\label{estI1}
&I_1 
= 	\varepsilon\sum_{k=0}^{2^N-1}{\E}^{\varepsilon, N}\Big  [\check Y_{\frac{t_k}{\varepsilon}}^{N,k} - \check Y^{N,k}_\frac{t_{k+1}}{\varepsilon}\Big ]\\
&\leq \varepsilon c\sum_{k=0}^{2^N-1}{\E}^{\varepsilon, N}\Big  [\Big(1+ |\tilde Z_{t_k}^N| + \Big(\int_{H \setminus 0}|\tilde U_{t_k}^N(\cdot)|^2 \nu_1(d\ell) \Big)^\frac{1}{2}\Big) (1+|\hat {\mathcal{Q}}_{\frac{t_k}{\varepsilon}}^{N}|_K+|\hat {\mathcal{Q}}_{\frac{t_k+1}{\varepsilon}}^{N}|_K)\Big ],\\
&\leq \varepsilon c\sum_{k=0}^{2^N-1}\Big({\E}^{\varepsilon, N}\Big  [  \rho\Big(1+ |\tilde Z_{t_k}^N|^2 + \Big(\int_{H \setminus 0}|\tilde U_{t_k}^N(\cdot)|^2 \nu_1(d\ell) \Big) \Big ]\Big)^{1/2}\cdot\\
&\cdot \Big({\E}^{\varepsilon, N}\Big  [  \rho^{-1}\Big (1+|\hat {\mathcal{Q}}_{\frac{t_k}{\varepsilon}}^{N}|_K+|\hat {\mathcal{Q}}_{\frac{t_k+1}{\varepsilon}}^{N}|_K\Big )^2\Big ]\Big)^{1/2}, \\
& \leq \varepsilon c \,  \kappa_1^\frac{1}{2}\sum_{k=0}^{2^N-1}\Big({\E}^{\varepsilon, N}\Big  [  \rho\Big(1+ |\tilde Z_{t_k}^N|^2 + \Big(\int_{H\setminus 0}|\tilde U_{t_k}^N(\cdot)|^2 \nu_1(d\ell) \Big) \Big ]\Big)^{1/2}\cdot\\
& \cdot \Big({\E}^{\varepsilon, N}\Big  [  \tilde \rho_1(1)\Big (1+|\hat {\mathcal{Q}}_{\frac{t_k}{\varepsilon}}^{N}|_K+|\hat {\mathcal{Q}}_{\frac{t_k+1}{\varepsilon}}^{N}|_K\Big )^2\Big ]\Big)^{1/2},
\end{align}
where $  \tilde \rho_1(1)$ is defined in \eqref{expmart} and $\kappa_p$ is the constant defined in \eqref{stima_path} with $p=1$,
and $c>0$ is the constant (independent of $k$ and $N$) appearing in \eqref{star1_2}.

We define a new probability $\tilde{\mathbb{P}}^{\varepsilon ,N}$ by $d \tilde{\mathbb{P}}^{\varepsilon ,N}=\tilde{\rho}_1(1) d \, \mathbb{P}^{\varepsilon ,N}$. Then last term in the previous inequality can be written as:
\begin{equation}
{\E}^{\varepsilon, N}\Big  [  \tilde \rho_1(1)\Big (1+|\hat {\mathcal{Q}}_{\frac{t_k}{\varepsilon}}^{N}|_K+|\hat {\mathcal{Q}}_{\frac{t_k+1}{\varepsilon}}^{N}|_K\Big )^2\Big ]=
\tilde{\E}^{\varepsilon, N}\Big  [ \Big (1+|\hat {\mathcal{Q}}_{\frac{t_k}{\varepsilon}}^{N}|_K+|\hat {\mathcal{Q}}_{\frac{t_k+1}{\varepsilon}}^{N}|_K\Big )^2\Big ]
\end{equation}
Under $\tilde \P^{\varepsilon,N}$,  
	$$
	\check{W}^1_s= \tilde W^1_s + \int_0^s \delta^{1,\varepsilon}_r \, dr, \quad s \in [0,1], 
	$$
	is a cylindrical Wiener process, and the compensator of $N_1(\omega, dt\,d\ell)$ is $(\gamma^{1,\varepsilon,N}_t( \omega, \ell) +1)^2\nu_1(d\ell) dt$ on $[0,1]$, while $\tilde{W}_2$ remains a Wiener process and $ N_2^\varepsilon(\omega, dt\,d\ell)$ has compensator $\displaystyle(\gamma^{2,\varepsilon,N}_\frac{t}{\varepsilon} (\omega, \ell) +1) \frac{\nu_2(d\ell)}{\varepsilon}\, dt$.
 
Therefore the forward equation \eqref{F-B_limit_system_3} reads:
\begin{align}\label{forward noise nuovo}
\left\{
\begin{array}{llll}
&d X_t=   A X_t dt - 2 R \delta^{1, \varepsilon}_{\frac{t}{\varepsilon}}dt
+  \int_{H \setminus 0} w(2 \,  \gamma_{\frac{t} {\varepsilon}}^{1, \varepsilon}(w) +  \gamma_{\frac{t} {\varepsilon}}^{1, \varepsilon}(w)^2)\nu_1(dw) dt +
  R d \check{W}_t^1 + \\ &\int_{H \setminus 0} w [N_1(dt\,dw) - ( \gamma_{\frac{t}{\varepsilon}}^{1, \varepsilon}(w)+1)^2
\nu_1(dw) dt], \\
&X_0 = x_0,\\
\end{array}
\right.
\end{align}
moreover we have that, by \eqref{mathcalQ^N} and \eqref{eqfast_changedtime_2}, 
%\begin{align*}
%		\hat {\mathcal{Q}}_s^N = \sum_{k=0}^{2^N-1}\hat {\mathcal{Q}}_s^{N,k} \,1_{\big[\frac{t_{k}}{\varepsilon},\frac{t_{k+1}}{\varepsilon}\big ]}(s),\quad s \in \Big[0, \frac{1}{\varepsilon}\Big ], 
%	\end{align*}
 $(\hat {\mathcal{Q}}_t^{N})_{t \in [0, 1/\varepsilon]}$ satisfies
 % , namely for $k=0,..., 2^N-1$,
\begin{align*}
\left\{
\begin{array}{ll}
d \hat {\mathcal{Q}}_t^{N}&=   [B  \hat {\mathcal{Q}}^{N}_t + F(X_{\varepsilon t}, \hat {\mathcal{Q}}_t^{N})]  dt -  G \delta^{2, \varepsilon, N}_t\, dt\!\!-\!\int_{K \setminus 0} w  \, \gamma_{t}^{2, \varepsilon}(w)
\nu_2(dw) dt+ G  d\hat {\tilde W}^2_t\\
&  +\int_{K\setminus 0} w [N_2(dt\,dw) - (1+ \gamma_t^{2,\varepsilon}(\omega))\nu_2(dw) dt]
\\
\hat {\mathcal{Q}}^{N}_0 &= q_0,
\end{array}
\right.
\end{align*}
with  
$\hat {\tilde W}^2_t = \frac{1}{\sqrt \varepsilon} d\tilde W_{\varepsilon t}^2$. On the other hand, according to \eqref{fbsystem_2},  $(\hat Q_t^{\varepsilon} = Q^{\varepsilon}_{\varepsilon t})_{ t \in [0,1/\varepsilon]}$ satisfies
\begin{align*}
\left\{
\begin{array}{llllll}
d \hat Q^\varepsilon_t& = [B \hat Q_t^\varepsilon + F(X_{\varepsilon t}, \hat Q_t^\varepsilon)]dt - G \delta^{2, \varepsilon, N}_t\, dt\!\!-\!\int_{K \setminus 0} w  \, \gamma_{t}^{2, \varepsilon}(w)
\nu_2(dw) dt+   G  d\hat {\tilde W}^2_t\\
& +\int_{K\setminus 0} w [N_2(dt\,dw) - (1+ \gamma_t^{2,\varepsilon}(\omega))\nu_2(dw) dt]
\\
\hat Q_0^\varepsilon& =q_0.
\end{array}
\right.
\end{align*}
Then by Lemma \ref{L:fasteq} with $g_t = - G \delta^{2, \varepsilon, N}_t\!-\!\int_{K \setminus 0} w  \, \gamma_{t}^{2, \varepsilon}(w)
\nu_2(dw) $, and recalling that \eqref{unifbounddelta}, 
there exists $k_2 < \infty$, independent of $\varepsilon$ and $N$  such that,
\begin{align}\label{supestQBIS}
	\sup_{t \in [0, 1/\varepsilon]}\tilde{\E}^{\varepsilon, N}[|\hat {\mathcal Q}^N_t|_K^2] 
	&\leq k_2 \Big(1+ |q_0|_K^2  + \sup_{t \in [0, 1/\varepsilon]}\tilde{\E}^{\varepsilon, N}[|X_{\varepsilon t}|_H^2]+|G|^2_{L(\Xi, K)}\Big).
\end{align}

%Arguing as before, it follows that, for all $p \geq 1$,  \textcolor{red}{non capisco più il senso di questo conto, secondo me non serve il cambio di probabilità che non agisce sul noise della veloce}
%\begin{align*}
%\E^{\varepsilon, N}\Big[\Big|\int_0^t e^{(t-s)B} G d \hat {\tilde L}_s^2\Big |_K^p\Big ]&= \E^{\varepsilon, N}\Big[\rho^{-1/2}\rho^{1/2}\Big|\int_0^t e^{(t-s)B} G d \hat {\tilde L}_s^2\Big |_K^p\Big ]\\
%&\leq 	(\E^{\varepsilon, N}[\rho^{-1}])^{1/2}\Big(\E^{\varepsilon, N}\Big[\rho\Big|\int_0^t e^{(t-s)B} G d \hat {\tilde L}_s^2\Big |_K^{2p}\Big ]\Big)^{1/2}\\
%&\leq \kappa_1^{1/2} \Big(\E\Big[\Big|\int_0^t e^{(t-s)B} G d \hat {L}_s^2\Big |_K^{2p}\Big ]\Big)^{1/2},
%\end{align*}
%where   $d\hat L^2_t = d\hat W^2_t + \int_{\Xi\setminus 0} w [N_2(dt\,dw) - \nu_2(dw) dt]$ with
%\hat W^2_t = \frac{1}{\sqrt \varepsilon} dW_{\varepsilon t}^2$.
Recalling estimate  \eqref{crescitaLpX}, 
\begin{align*}
	\sup_{ t \in [0,\frac{1}{\varepsilon}]}\tilde{\E}^{\varepsilon, N}[|X_{\varepsilon t}|_H^2]&
	\leq 	c \Big( 1+ |x_0|^2]\Big).
\end{align*}
where the constant $c$, in view of \eqref{stime unif girsanov}, can be chosen independently of $\varepsilon$, $N$,
% Therefore, recalling that  by Hypothesis 
% \noindent {\sc{\textbf{(H$\beta^B$)}}} \textcolor{red}{DEVE VALERE PER UN PROCESSO CON COMPENSATORE $(1+ \gamma)$ }  \ref{R:Levyconv} 
% \begin{equation}\label{estbetapB}
% 	\sup_{t \geq 0} \E\Big[\Big|\int_0^t e^{(t-s)B} G d \hat  L_s^2\Big |_K^2\Big ] < \infty, 
% \end{equation}
estimate \eqref{supestQBIS} reads \begin{align}\label{supestQTRIS}
	\sup_{t \in [0, 1/\varepsilon]}\E^{\varepsilon, N}[|\hat {\mathcal Q}^N_t|_K^2] \leq c.
\end{align}
%for some constant $c >0$ independent of $\varepsilon$ and $N$. 

Then, plugging  \eqref{supestQTRIS} in 
\eqref{estI1} we get
\begin{align}\label{estI1BIS}
I_1 &\leq 3\varepsilon c\sum_{k=0}^{2^N-1}\Big({\E}^{\varepsilon, N}\Big  [ \Big(1+ |\tilde Z_{t_k}^N|^2 + \int_{H \setminus 0}|\tilde U_{t_k}^N(\cdot)|^2 \nu_1(d\ell) \Big) \Big ]\Big)^{1/2}. 
\end{align}
%We recall that the laws of $\tilde Z^N$ and $\tilde U^N(\cdot)$ depend only on the law of $X$, so in particular
%\begin{align}
%		\E^{\varepsilon, N}\Big[\rho \int_0^1 |\tilde Z^N_t|_{\Xi^\ast}^2 dt\Big]&=\E\Big[ \int_0^1 |\tilde  Z^N_t|_{\Xi^\ast}^2 dt\Big]< + \infty,\label{rhoZtilde} \\
%		\E^{\varepsilon, N}\Big[\rho \int_0^1 \int_{\Xi\setminus 0} |\tilde  U^N_t(w)|^2 \nu_1(dw) dt\Big]&=\E\Big[ \int_0^1 \int_{\Xi\setminus 0} |\tilde  U^N_t(w)|^2 \nu_1(dw) dt\Big]< + \infty.\label{rhoUtilde}
%	\end{align}
	Proceeding as above, recalling \eqref{rhoZ}-\eqref{rhoU} and \eqref{ZN}-\eqref{UN},  we get 
	\begin{align}
	&{\E}^{\varepsilon, N}\Big  [\rho\Big(1+ |\tilde Z_{t_k}^N|^2 + \int_{H \setminus 0}|\tilde U_{t_k}^N(\cdot)|^2 \nu_1(d\ell) \Big) \Big ]\notag\\
	&\leq c \Big({\E}\Big  [1+ 2^{N}\int_{t_{k-1}}^{t_{k}}|\bar Z_{t}|^2 dt + 2^{N}\int_{t_{k-1}}^{t_{k}}\int_{H \setminus 0}|\bar U_{t}^N(\cdot)|^2 \nu_1(d\ell)  \Big ]\Big ).\label{finalestI1}
	\end{align}
 with $c$ doesn't  depend neither on $\varepsilon$ nor on $N$.
	Plugging in turn \eqref{finalestI1} in \eqref{estI1BIS} we achieve
	\begin{align}\label{estI1TRIS}
I_1 &\leq \varepsilon C\sum_{k=0}^{2^N-1}\Big({\E}\Big  [1+ 2^{N}\int_{t_{k-1}}^{t_{k}}|\bar Z_{t}|^2 dt + 2^{N}\int_{t_{k-1}}^{t_{k}}\int_{H \setminus 0}|\bar U_{t}^N(\cdot)|^2 \nu_1(d\ell)  \Big ]\Big )^{1/2}\\
&\leq \varepsilon C\sum_{k=0}^{2^N-1}\Big(1+2^{\frac{N}{2}}\Big({\E}\Big  [\int_{t_{k-1}}^{t_{k}}|\bar Z_{t}|^2 dt \Big]\Big)^{1/2}+ 2^{\frac{N}{2}}\Big({\E}\Big  [\int_{t_{k-1}}^{t_{k}}\int_{H \setminus 0}|\bar U_{t}^N(\cdot)|^2 \nu_1(d\ell)  \Big ]\Big )^{1/2}\Big)\\
&\leq \varepsilon C \Big(2^N+2^{\frac{3}{2}N}\Big({\E}\Big  [\int_{0}^{1}|\bar Z_{t}|^2 dt \Big]\Big)^{1/2}+ 2^{\frac{3}{2}N}\Big({\E}\Big  [\int_{0}^{1}\int_{H \setminus 0}|\bar U_{t}^N(\cdot)|^2 \nu_1(d\ell)  \Big ]\Big )^{1/2}\Big)\\
&= C \varepsilon (2^N + 2^{\frac{3}{2}N}). 
\end{align}

%%%%%%%%%%%%%%%%%%%%%%%%%%%%%%%%%
 Concerning now $I_5$, we notice that $\P$-a.s., for all $s >0$,  
\begin{align}\label{diffestQBIS}
|\hat Q^\varepsilon_s - \hat {\mathcal Q}^N_{s}|\leq K\int_0^s e^{-\mu (s-t)}|X_{\varepsilon t}- X^N_{\varepsilon t}|\,dt,
\end{align}
where $\mu$ is the dissipative constant in {\sc{\textbf{(HF+B)}}}, and $K$ does not depend of $\varepsilon$ and $N$.
Thus we deduce, thanks to \eqref{diffestQBIS} and proceeding as for $I_4$, 
\begin{align}\label{estI6BIS}
I_5 &\leq  L_q K  {\E}^{\varepsilon, N} \Big [\int_0^{1}\Big (1 + |\bar Z_{t}|^2_{\Xi^\ast}+\int_{H \setminus 0}|\bar U_{t}(w)|^2 \nu_1(dw)\Big )\Big(\int_0^{\frac{t}{\varepsilon}} e^{-\mu (s-r)}|X_{\varepsilon r}- X^N_{\varepsilon r}|\,dr \Big) dt\Big ]\notag\\
&\leq   C  ({\E}^{\varepsilon, N} [\Delta_{X, N}^4])^{1/4}. 
\end{align}

%%%%%%%%%%%%%%%%%%%%%%%%%%%

Collecting \eqref{estI4}-\eqref{I2est}-\eqref{estI3}-\eqref{estI6BIS}-\eqref{estI1TRIS},  we get from \eqref{lastid} that
\begin{align}
	&Y^\varepsilon_0 -\bar Y_0\leq  C(\varepsilon (2^N + 2^{\frac{3}{2}N}) + ({\E} [\Delta_{Z, N}])^{1/2}+ ({\E} [\Delta_{U, N}])^{1/2} +  ({\E}[\Delta_{X, N}^4])^{1/4}).\label{lastidBIS}
\end{align}
 Finally we notice that {having $X$ càdlàg paths} and since by \eqref{crescitaLpX} $\E[\sup_{t \in [0,1]}|X_t|^4]< \infty$, 
 $$
 {\E}[\Delta_{X, N}^4]\rightarrow 0 \,\,\textup{as}\,\,N \rightarrow \infty. 
 $$
 Moreover, by \eqref{convZN}-\eqref{convUN}, 
  $$
 {\E}[\Delta_{Z, N}]\rightarrow 0 \quad \textup{and} \quad  {\E}[\Delta_{U, N}]\rightarrow 0\,\,\textup{as}\,\,N \rightarrow \infty. 
 $$
 Therefore, letting first $\varepsilon\rightarrow 0$  and then $N\rightarrow \infty$  in \eqref{lastidBIS} we get the result. 
\endproof

\section{The two scale control problem}
One of the purposes of the article is to give a representation of the limit of the value functions of a sequence of optimal control problems for a singularly perturbed infinite dimensional state equation.

\subsection{Formulation of the problem}
We denote by  $\mathcal A$  the set of predictable processes $(\alpha_t)_{t \in [0,1]}$ taking values in a complete metric space $U$. 
Given a solution $(X, Q^\varepsilon)$ of \eqref{SDE} and 
a control $\alpha \in \mathcal A$,   we associate to them a probability measure $\P^{\varepsilon, \alpha}$ on $(\Omega, \mathcal F)$ that  we are going to describe. 
Let  $b:H \times K\times U\rightarrow H$, $\rho: U \rightarrow K$, $r:H \times K\times U \times H\rightarrow \R$ and $\gamma: U\times K\rightarrow \R$ be measurable functions satisfying the following assumptions. 

\medskip

\noindent {\sc{\textbf{(H$b\rho$)}}} There exist positive constants   $L_b$ and $M'$ such that, for every $x, x' \in H$, $q, q' \in K$, $a \in U$,  
\begin{align*}
	|b(x,q,a) - b(x', q',u)|_H &\leq L_b (|x-x'|_H+|q-q'|_K),\\
	|b(x, q, a)|_H + |\rho(a)|_K & \leq M'. 
\end{align*}

\medskip

\noindent {\sc{\textbf{(H$r\gamma$)}}} There exist  constants   $C_r >1$, $L_r >0$ and $C_\gamma >1$ such that, for every $x\in H$, $q\in K$, $a \in U$, $ w \in \Xi$,  
\begin{align*}
	0<\eta\leq  r(x,q,a,w) &\leq C_r, \\
	|r(x,q,a,w)-r(x',q',a,w)| &\leq L_r(|x-x'|_H+|q-q'|_K),\\
	0<\eta'\leq \gamma(a,w) &\leq C_\gamma. 
\end{align*}

We also need to add the following additional assumption on the operator $R$. 

\medskip

\noindent {\sc{\textbf{(H')}}} $R$ admits a bounded right inverse $R^{-1} \in L( H,\Xi)$. 

\medskip

%Let $A$, $B$, $R$,   $G$, and $F$ be   given as in Section \ref{S:3},  and
Consider the  processes  
\begin{align*}
	M^1&:=\int_{[0,\,\cdot]}R^{-1} b(X_t, Q^\varepsilon_t, \alpha_t) d W^1_t + \int_{]0,\,\cdot]} \int_{H \setminus 0}[{r(X_t, Q^\varepsilon_t, \alpha_t,w)}-1] (N_1(dt\,dw)-\nu_1(dw)dt),\\ 
	M^2&:=\int_{[0,\,\cdot]}\frac{\rho(\alpha_t)}{\sqrt \varepsilon  } d W^2_t + \int_{]0,\,\cdot]} \int_{K \setminus 0}[{\gamma(\alpha_t,w)}-1] \Big (N^{\varepsilon}_2(dt\,dw)-\frac{\nu_2(dw)}{\varepsilon}dt\Big ),
\end{align*}
and set  
	$$
	\frac{d \P^{\varepsilon, \alpha}}{d\P}:= \mathcal E_1(M^1+M^2).
	$$
By the Girsanov Theorem \ref{Girsanov}, 
the compensator of $(N_1(dt\,dz), N^\varepsilon_2(dt\,dz))$ under $\P^{\varepsilon, \alpha}$ is $({r(X_t, Q^\varepsilon_t, \alpha_t,w)}\nu_1(dw) dt, \frac{{\gamma(\alpha_t,w)}}{\varepsilon}\nu_2(dw) dt)$, and the process
	\begin{align*}
	(\tilde W^1,\tilde W^2)&:=\left(W^1 - \int_{[0,\,\cdot]}R^{-1} b(X_t, Q^\varepsilon_t, \alpha_t) dt, W^2 - \int_{[0,\,\cdot]}\frac{\rho(\alpha_t)}{\sqrt \varepsilon}   dt\right)
	%(N_1^{\P^{\varepsilon, \alpha}}(dt\,dz),N_2^{\varepsilon,\P^{\varepsilon, \alpha}}(dt\,dz))&:=,
\end{align*}
is such that $\tilde W^1$ and $\tilde W^2$ are independent cylindrical Wiener processes under $\P^{\varepsilon, \alpha}$.  
We can therefore conclude that the solution $(X, Q^\varepsilon)$ of \eqref{SDE} satisfies  as well the system of controlled SDEs: \begin{align}\label{ControlledSDE}
\left\{
\begin{array}{llll}
d X_t=   A X_t dt  + b(X_t, Q^\varepsilon_t, \alpha_t) dt +R d \tilde W^1_t + \int_{H \setminus 0} w [N_1(dt\,dw) - {r(X_t, Q^\varepsilon_t, \alpha_t,w)}\nu_1(dw) dt], \\
X_0 = x_0 \in H,\\
\varepsilon dQ^\varepsilon_t = [B Q_t^\varepsilon + F(X_t, Q_t^\varepsilon)]dt + G \rho(\alpha_t)dt+ \sqrt \varepsilon G d \tilde W^2_t\\
\qquad \qquad + \int_{K\setminus 0} w [\varepsilon N_2^\varepsilon(dt\,dw) - {\gamma(\alpha_t,w)}\nu_2(dw) dt],\\
Q_0^\varepsilon =q_0 \in K. 
\end{array}
\right.
\end{align}
\subsection{Solution of the problem: the BSDE approach}
Let us now consider a measurable function $l: H \times K \times U \rightarrow \R $   satisfying the following assumption. 

\medskip

\noindent {\sc{\textbf{(H$l$)}}} There exist positive constants   $L_l$ and $M_l$ such that, for every $x, x' \in H$, $q, q' \in K$, $a \in U$,  
\begin{align*}
	|l(x,q,a) - l(x', q',a)| &\leq L_l (|x-x'|_H+|q-q'|_K),\\
	|l(x, q, a)|  & \leq M_l. 
\end{align*}

%\noindent {\sc{\textbf{(H$h$)}}} There exists positive constants   $L_h$ and $M_h$ such that, for every $x, x' \in H$, $q, q' \in K$, $a \in U$,  
%\begin{align*}
%	|h(x) - h(x')| &\leq L_h (|x-x'|_H),\\
%	 |h(x)| & \leq M_h. 
%\end{align*}
%
%\medskip 

We define, for $x \in H$, $q \in K$, $ z, \zeta \in \Xi^\ast$, $u \in L^2(\nu_1)$, $\theta \in L^2(\nu_2)$, the function 
\begin{align}\label{hamiltonian}
\psi(x,q,z,\zeta, u, \theta)&:= \inf_{a \in U} \Big \{l(x,q, a) + z[R^{-1}b(x,q,a)]+ \zeta \rho(a) \\
&+ \int_{H \setminus 0}u(w)(r(x,q,a,w)-1) \nu_1(dw)+ \int_{K \setminus 0}\theta(w)(\gamma(a,w)-1) \nu_2(dw) \Big \}. 
\end{align}
\begin{proposition}\label{P:f_psi}
		Let assumptions {\sc{\textbf{(H$b\rho$)}}},  {\sc{\textbf{(H$r\gamma$)}}}, {\sc{\textbf{(H')}}} and {\sc{\textbf{(H$l$)}}}  hold true. Then the function $\psi$ in \eqref{hamiltonian} verifies hypothesis {\sc{\textbf{(H$\psi$)}}}.
\end{proposition}
\proof 
% First of all,  $\psi$ is measurable, and
% 		$
% 		\sup_{x \in H, q \in K}\psi(x, q, 0,0,0,0) %= \inf_{a \in U} l(x,q, a)  
% 		\leq M.
% 		$
% 		Moreover, 
% 		\begin{align*}
% 	&		l(x,q, a) + z[R^{-1}b(x,q,a)]+ \zeta \rho(a) \\
% &+ \int_{\Xi \setminus 0}u(w)(r(x,q,a,w)-1) \nu_1(dw)+ \int_{\Xi \setminus 0}\theta(w)(\gamma(a,w)-1) \nu_2(dw)\\
% &=l(x',q', a) + z'[R^{-1}b(x',q',a)]- \zeta' \rho(a) \\
% &+ \int_{\Xi \setminus 0}u'(w)(r(x',q',a,w)-1) \nu_1(dw)+ \int_{\Xi \setminus 0}\theta'(w)(\gamma(a,w)-1) \nu_2(dw)\\
% &+ (l(x,q, a) -l(x',q', a) ) + (z-z')[R^{-1}b(x,q,a)]+ z'[R^{-1}(b(x,q,a)-b(x',q',a))]\\
% &+ \int_{\Xi \setminus 0}(u(w)-u'(w))(r(x,q,a,w)-1)\nu_1(dw)+\int_{\Xi \setminus 0}u'(w)(r(x,q,a,w)-r(x',q',a,w))\nu_1(dw)\\
% %&+\int_{\Xi \setminus 0}(u'(w)-u(w)) \nu_1(dw)	
% &+ \int_{\Xi \setminus 0}(\theta(w)-\theta'(w))(\gamma(a,w)-1) \nu_2(dw)\\
% &\leq l(x',q', a) + z'[R^{-1}b(x',q',a)]- \zeta' \rho(a) \\
% &+ \int_{\Xi \setminus 0}u'(w)(r(x',q',a,w)-1) \nu_1(dw)+ \int_{\Xi \setminus 0}\theta'(w)(\gamma(a,w)-1) \nu_2(dw)\\
% &+L_l (|x-x'|_H+|q-q'|_K) + ||R^{-1}||_{L(H, \Sigma)}M'|z-z'|_{\Xi^\ast}\\
% &+||R^{-1}||_{L(H, \Sigma)}L_b |z'|_{\Xi^\ast}(|x-x'|_H+|q-q'|_K)\\
% &+ (C_r+1) \sqrt{\nu_1(\Sigma)}\Big(\int_{\Xi \setminus 0}|u_s(w)- u'_s(w)|^2 \nu_1(dw)\Big)^\frac{1}{2}\\
% &+ (C_\gamma+1)\sqrt{\nu_1(\Sigma)}\Big(\int_{\Xi \setminus 0}|\theta_s(w)- \theta'_s(w)|^2 \nu_2(dw)\Big)^\frac{1}{2}\\
% &+ L_r(|x-x'|_H+|q-q'|_K)\sqrt{\nu_1(\Sigma)}\Big(\int_{\Xi \setminus 0}|u'_s(w)|^2 \nu_1(dw)\Big)^\frac{1}{2}. 
% 	\end{align*}
By standard estimates 
\begin{align*}
	&		|\psi(x,q,z,\zeta,u, \theta)-\psi(x',q',z',\zeta',u', \theta')| \\
&\leq  L_u \Big(\int_{\Xi \setminus 0}|u_s(w)- u'_s(w)|^2 \nu_1(dw)\Big)^\frac{1}{2} + L_\theta \Big(\int_{\Xi \setminus 0}|\theta_s(w)- \theta'_s(w)|^2 \nu_2(dw)\Big)^\frac{1}{2}\\
		&+ L_x \Big(1 + |z|_{\Xi^\ast} + \Big(\int_{H \setminus 0}|u_s(w)|^2 \nu_1(dw)\Big)^\frac{1}{2}\Big)|x-x'|_H\\
		&+L_q \Big(1 + |z|_{\Xi^\ast} + \Big(\int_{H\setminus 0}|u_s(w)|^2 \nu_1(dw)\Big)^\frac{1}{2}\Big)|q-q'|_K, 
	\end{align*}
	for some constants $L_x, L_q, L_z, L_\zeta,  L_u, L_\theta>0$ defined opportunely.
	
	Finally, we have to verify that 
for every $x \in H$, $q \in K$, $z, \zeta \in \Xi^\ast$, $u, u'\in L^2(\nu_1)$, $\theta, \theta' \in L^2(\nu_2)$, there exist measurable functions $\gamma_1: H \rightarrow \R$, $\gamma_2: K \rightarrow \R$  (depending on $x$, $q$, $z,\zeta$, $u$, $\theta, \theta'$) such that 
	\begin{align}
	\psi(x,q,z,\zeta,u, \theta)-\psi(x,q,z,\zeta,u', \theta) \leq \int_{H \setminus 0} (u(w) -u'(w)) \gamma_1(w) \nu_1(dw)\label{1toprove}\\
		\psi(x,q,z,\zeta,u, \theta)-\psi(x,q,z,\zeta,u, \theta') \leq \int_{K\setminus 0} (\theta(w) -\theta'(w)) \gamma_2(w) \nu_2(dw),\label{2toprove}
	\end{align}
	and satisfying  
	\begin{align*}
	&C_1 (1 \wedge |w|_H) \leq \gamma_1(w) \leq C_2 (1 \wedge |w|_H), \quad C_1 \in (-1, 0], C_2 \geq 0, \\
	&\bar C_1 (1 \wedge |w|_K) \leq \gamma_2(w) \leq \bar C_2 (1 \wedge |w|_K), \quad \bar C_1 \in (-1, 0], \bar C_2 \geq 0.
	\end{align*}
	We have
	\begin{align}\label{hamiltonian_diff_u}
&\int_{H \setminus 0}u(w)(r(x,q,a,w)-1) \nu_1(dw)\\
&= \int_{H \setminus 0}u'(w)(r(x,q,a,w)-1) \nu_1(dw)+\int_{H \setminus 0}(u(w)-u'(w))(r(x,q,a,w)-1) \nu_1(dw)\\
&\leq  \int_{H \setminus 0}u'(w)(r(x,q,a,w)-1) \nu_1(dw)+\int_{H \setminus 0}(u(w)-u'(w))\gamma_1(w) \nu_1(dw)
%&\leq \int_{\Xi \setminus 0}u'(w)(r(x,q,a,w)-1) \nu_1(dw)+(C_r-1)\int_{\Xi \setminus 0}(u(w)-u'(w)) \nu_1(dw).
\end{align} 
with $\gamma_1(w) := \sup_{a \in U}(r(x,q,a,w)-1) I_{[0,+\infty[}(u(w)-u'(w)) + \inf_{a \in U}(r(x,q,a,w)-1) I_{]0,+\infty[}(u'(w)-u(w))$. Notice that, by assumption {\sc{\textbf{(H$r\gamma$)}}},  
$$
-1<\eta-1\leq \gamma_1(w)\leq C_r-1, \quad C_r >1. 
$$
Adding and subtracting in both sides $l(x,q, a) + z[R^{-1}b(x,q,a)]+ \zeta \rho(a) + \int_{H \setminus 0}u(w)(r(x,q,a,w)-1) \nu_1(dw)$ and taking the infimum  over $a \in U$, we get \eqref{1toprove}. 

Analogously, one observes that 
\begin{align*}
	&\int_{K\setminus 0}\theta(w)(\gamma(a,w)-1) \nu_2(dw) \leq \int_{K \setminus 0}\theta'(w)(\gamma(a,w)-1) \nu_2(dw)\\
	&+ \int_{K \setminus 0}(\theta(w)-\theta'(w))\gamma_2(w) \nu_2(dw)
\end{align*}
with $\gamma_2(w) := \sup_{a \in U}(\gamma(a,w)-1)  I_{[0,+\infty[}(\theta(w)-\theta'(w)) +\inf_{a \in U}(\gamma(a,w)-1)  I_{]0,+\infty[}(\theta'(w)-\theta(w)) ]$. In particular, by assumption {\sc{\textbf{(H$r\gamma$)}}},  
$$
-1<\eta'-1\leq \gamma_2(w)\leq C_\gamma-1, \quad C_\gamma >1. 
$$
 Proceeding as before \eqref{2toprove} follows. 

\endproof

Let  $h:H \rightarrow \R $ be a measurable function satisfying {\sc{\textbf{(H$h$)}}}.
%Let $h: H \rightarrow \R$ be a Lipschitz continuous function with constant $L_h>0$. 
We consider the following BSDE: 
$\P$-a.s., for $t \in [0,1]$, 
 \begin{align}\label{BSDE_control}
	 Y_t^{\varepsilon} &=  h(X_1) + \int_t^1 \psi\Big (X_s, Q_s^{\varepsilon},Z^{\varepsilon}_s,  \frac{\Xi^{\varepsilon}_s}{\sqrt{\varepsilon}}, U^{\varepsilon}_s(\cdot), \frac{\Theta^{\varepsilon}_s(\cdot)}{\varepsilon}\Big ) ds - \int_t^1 Z^{\varepsilon}_s d W^1_s - \int_t^1 \Xi^{\varepsilon}_s d  W^2_s \\
	&- \int_t^1 \int_{H \setminus 0} U^{\varepsilon}_s(w) (N_1(ds \, dw) - \nu_1(dw) ds)- \int_t^1 \int_{K \setminus 0} \Theta^{\varepsilon}_s(w) \Big(N_2^\varepsilon(ds \, dw) - \frac{\nu_2(dw)}{\varepsilon} ds\Big).
	\end{align}

	\begin{lemma}\label{L:Vesp}
Let assumptions  {\sc{\textbf{(HAB)}}}, {\sc{\textbf{(HRG)}}}, {\sc{\textbf{(HF)}}}, {\sc{\textbf{(HF+B)}}}, {\sc{\textbf{(H$\beta^B$)}}},  {\sc{\textbf{(H$b\rho$)}}},  {\sc{\textbf{(H$r\gamma$)}}}, {\sc{\textbf{(H')}}},  {\sc{\textbf{(H$l$)}}} and  {\sc{\textbf{(H$h$)}}} hold true.
Then for every $\varepsilon >0$,  there exists a unique solution $(X, Q^\varepsilon, Y^\varepsilon, Z^\varepsilon, \Xi^\varepsilon, U^\varepsilon(\cdot))$   to the forward-backward system \eqref{BSDE_control}, and 
$$
		Y_0^\varepsilon=V^\varepsilon(x_0, q_0), \quad x_0 \in H, q_0 \in K, 
		$$
		where 
 $V^\varepsilon$ is the value function in  \eqref{valuefunc}. 
 	\end{lemma}
	\proof
The well-posedness of the forward-backward system \eqref{BSDE_control} directly follows from Theorem \ref{T:welpos1BIS} together with Lemma \ref{P:f_psi}. 	Then the control interpretation of $Y_0^\varepsilon$ follows as in \cite{bandiniconfortolacosso}, that  extends to the infinite dimensional framework the classical representation in finite dimension, see e.g. \cite{BaBuPa}.
	\endproof
	The following result is a direct consequence of Theorem \ref{T:convYesp}, together with Lemma \ref{L:Vesp}.
	\begin{theorem}\label{th:main_contr}
Let assumptions  {\sc{\textbf{(HAB)}}}, {\sc{\textbf{(HRG)}}}, {\sc{\textbf{(HF)}}}, {\sc{\textbf{(HF+B)}}}, {\sc{\textbf{(H$\beta^B$)}}},  {\sc{\textbf{(H$b\rho$)}}},  {\sc{\textbf{(H$r\gamma$)}}}, {\sc{\textbf{(H')}}},  {\sc{\textbf{(H$l$)}}} and  {\sc{\textbf{(H$h$)}}} hold true, and let $(X, \bar Y, \bar Z, \bar U(\cdot))$ be  the unique solution to the forward-backward system \eqref{F-B_limit_system_2}. Then 
$$
\lim_{\varepsilon \rightarrow 0}V^\varepsilon(x_0, q_0)= \bar Y_0, \quad x_0 \in H, q_0 \in K, 
$$
where $V^\varepsilon$, $\varepsilon >0$,  is the value function in  \eqref{valuefunc}.  
	\end{theorem}
	\section{Control interpretation of the limit forward-backward system}
Let $(X, \bar Y, \bar Z, \bar U(\cdot))$   be the unique solution to the forward-backward system \eqref{F-B_limit_system_2}. We aim now at interpreting $\bar Y_0$ as the value function of a corresponding reduced control problem. 
To this end, we start by noticing that, for every $t \in [0,1]$,  
	\begin{align}\label{BSDEpv}
\bar Y_t = h(X_1)+   \int_t^1\lambda(X_s,\bar Z_s,\bar U_s(\cdot)) ds 
	 - \int_t^1 \bar Z_s d  W^1_s - \int_t^1 \int_{H \setminus 0} \bar U_s(w) [N_1(ds \, dw) - \nu_1(dw) ds], 
\end{align}
where $(X_t)_{t \geq 0}$ is the mild solution of
\begin{align*}
\left\{
\begin{array}{llll}
d X_t=   A X_t dt + R d W_t^1 + \int_{H \setminus 0} w [N_1(dt\,dw) - \nu_1(dw) dt], \\
X_0 = x_0,
\end{array}
\right.
\end{align*}
and 
(see Theorem \ref{T:ergodic})
$\lambda : H  \times \Xi^\ast \times L^2(\nu_1) \rightarrow \R$ is uniquely determined and 
satisfies,  
 for fixed $x \in H$, $z, z' \in \Xi^\ast$, $u, u' \in L^2(\nu_1)$, 
\begin{align}
	|\lambda(x,z,u)-\lambda(x,z',u')|&\leq L_z |z-z'|_{\Xi^\ast}+ L_u\Big(\int_{H\setminus 0}|u(w)- u'(w)|^2 \nu_1(dw)\Big)^\frac{1}{2},\label{Lip_lambda}\\
	|\lambda(x,0,0)|&\leq c(1+ |x|_H),\label{bound_lambda}
\end{align}
for some positive constants  $L_z$,  $L_u$ and $c$.
Moreover,  {$\lambda(\cdot, z,u)$  is concave with respect to $z$ and $u$}, so by the Fenchel-Moreau theorem one can write $\lambda = \lambda_{\ast\ast}$, where, for all $x \in H$,
\begin{equation}\label{def_lambdastar}
\lambda_\ast(x,p,v)= \inf_{z \in \Xi^\ast, u \in L^2(\nu_1)}(-zp - \langle u, 1-v\rangle_{L^2(\nu_1)}-  \lambda(x,z,u)), \quad p \in \Xi, v \in L^2(\nu_1). 
\end{equation} 
It follows that, for all $x \in H$, $z \in \Xi$ and $u \in L^2(\nu_1)$, 
\begin{equation*}
	\lambda(x,z,u)= \inf_{p \in \Xi, \, v \in L^2(\nu_1)}(-zp - \langle u, 1-v\rangle_{L^2(\nu_1)}-  \lambda_\ast(x,p,v)). 
\end{equation*}

Moreover by the Lipschitzianity of $\lambda$ we can restrict the infimum to a bounded set, namely the following  holds. 
\begin{proposition}\label{P:74}
We have 
\begin{equation}\label{gooddef}
	\lambda(x,z,u)= \inf_{\underset{|p|\leq L_z, \, |1-v|_{L^2(\nu_1)}\leq L_u}{p \in \Xi, \, v \in L^2(\nu_1):}}(-zp - \langle u, 1-v\rangle_{L^2(\nu_1)}-  \lambda_\ast(x,p,v)),  
\end{equation}	
where $L_z$ and $L_u$ are the Lipschitz constants in \eqref{Lip_lambda}. 
\end{proposition}
\proof 
By \eqref{def_lambdastar} together with \eqref{bound_lambda} we get 
$$
\lambda_\ast(x,p,v) \leq -\lambda(x,0,0) \leq c (1+|x|), \quad p \in \Xi, v \in L^2(\nu_1). 
$$
Then, again by \eqref{def_lambdastar} together with \eqref{Lip_lambda}, we obtain that, for all $p \in \Xi, v \in L^2(\nu_1)$, 
\begin{align*}
\lambda_\ast(x,p,v)&\leq  \inf_{z \in \Xi^\ast, u \in L^2(\nu_1)}(-zp - \langle u, 1-v\rangle_{L^2(\nu_1)} + L_z |z|_{\Xi^\ast} + L_u |u|_{L^2(\nu_1)})+ c.
\end{align*}
It follows that $\lambda_\ast(x,p,v)= - \infty$ if $|p|_{\Xi}> L_z$ or $|1-v|_{L^2(\nu_1)}> L_u$, which implies \eqref{gooddef}. 
\endproof

We  now define 
\begin{align*}
\mathcal S:=\{(p_s,v_s(\cdot))_{s \in [0,1]}: &\, p \textup{ progressively measurable s.t. } |p_s|_{\Xi} \leq L_z, \\
&v(\cdot) \textup{ predictable s.t. } |1-v_s|_{L^2(\nu_1)} \leq L_u   \}. 
\end{align*}
For every  $(p,v(\cdot)) \in \mathcal S$, 
let us set 
\begin{equation}\label{Mpv}
	M:=-\int_{[0,\,\cdot]}p^*_t d W^1_t + \int_{]0,\,\cdot]} \int_{H \setminus 0}(v_t(w)-1) [N( dt\,dw)-\nu_1(dw)dt]. 
\end{equation} 
Being $(1-v_t(\omega,  \cdot))\in L^2(\nu_1)$ for all $(\omega, t) \in \Omega \times [0,\,1]$, and since, for all $t \in [0,1]$,  
$$
\E\Big[e^{\frac{1}{2} \langle M^c\rangle_t +   \langle M^d\rangle_t}\Big]= \E\Big[e^{\frac{1}{2} \int_0^t |p|^2_s ds +   \int_0^t \int_{H \setminus 0} (1-v_s(w))^2 \nu_1(dw) ds}\Big]\leq \E\Big[e^{(\frac{1}{2} L_z^2  +   L_u) t}\Big]< \infty,  
$$ 
by Theorem \ref{Girsanov_gen} together with Corollary \ref{C_1543}, 
under  $\P^{p,v}$ defined by 
	$$
	\frac{d \P^{p,v}}{d\P}:= \mathcal E_1\left(M\right),
	$$
	the process
$$
W^p_t := \int_0^t p_s ds + W_t^1, \quad t \in [0,1], 
$$
is a $Q_1$-Wiener process, and $N_1(dt\,dw)$  has compensator $v_t(w) \nu_1(dw) dt$.

We can finally give the following result that follows the classical argument relating the BDSEs to control problems. We sketch the proof for the reader convenience.
\begin{theorem}\label{th:repr_as_control}
Assume that   {\sc{\textbf{(HAB)}}}, {\sc{\textbf{(HRG)}}}, {\sc{\textbf{(HF)}}}, {\sc{\textbf{(HF+B)}}}, {\sc{\textbf{(H$\beta^B$)}}},  {\sc{\textbf{(H$b\rho$)}}},  {\sc{\textbf{(H$r\gamma$)}}}, {\sc{\textbf{(H')}}},  {\sc{\textbf{(H$l$)}}} and  {\sc{\textbf{(H$h$)}}} hold true. Then 
$$
\bar Y_0 = \inf_{(p,v(\cdot)) \in \mathcal S} \E^{p,v}\Big[h(X_1) - \int_0^1 \lambda_\ast(X_s, p_s, v_s(\cdot))ds\Big].
$$
\end{theorem}
%\begin{remark}
%$(X_t)_{t \geq 0}$ is the mild solution of the controlled SDE
%\begin{align*}
%\left\{
%\begin{array}{llll}
%	d X_t=   A X_t dt - R p_t dt+ R d W_t^p + R\int_{\Xi \setminus 0} w [N_1(dt\,dw) - v_t(w)\nu_1(dw) dt]\\
%	\qquad \quad + R\int_{\Xi \setminus 0} w (v_t(w)-1)\nu_1(dw) dt, \\
%X_0 = x_0.
%\end{array}
%\right.
%\end{align*}
%\end{remark}
\proof
We prove the two inequalities separately.
 
\noindent $(\Leftarrow)$ Let $(p,v(\cdot)) \in \mathcal S$. By Proposition \ref{P:74}, for all $t \in [0,1]$, from \eqref{BSDEpv} 
we get 
	\begin{align*}
\bar Y_t &= h(X_1)+   \int_t^1\lambda(X_s,\bar Z_s,\bar U_s(\cdot)) ds 
	 - \int_t^1 \bar Z_s d  W^1_s - \int_t^1 \int_{H \setminus 0} \bar U_s(w) [N_1(ds \, dw) - \nu_1(dw) ds]\\
	 &\leq  h(X_1)-   \int_t^1\lambda_\ast(X_s,p_s,v_s(\cdot)) ds - \int_t^1 \bar Z_s p_s ds- \int_t^1 \int_{H \setminus 0} (1-v_s(w))\bar U_s(w) \nu_1(dw) ds
	 \\
	 &- \int_t^1 \bar Z_s d  W^1_s - \int_t^1 \int_{H \setminus 0} \bar U_s(w) [N_1(ds \, dw) - \nu_1(dw) ds]\\
	 &= h(X_1)-   \int_t^1\lambda_\ast(X_s,p_s,v_s(\cdot)) ds 	- \int_t^1 \bar Z_s d  W^p_s \\
	 & - \int_t^1 \int_{H \setminus 0} \bar U_s(w) [N_1(ds \, dw) - v_s(w)\nu_1(dw) ds].
\end{align*}
Applying the conditional expectation with respect to $\P^{p,v}$ we get 
	\begin{align*}
\bar Y_t 
	 &\leq  \E^{p,v}\Big [h(X_1)-   \int_t^1\lambda_\ast(X_s,p_s,v_s(\cdot)) ds\Big|\mathcal F_t \Big], \quad t \in [0,1]. 
	 \end{align*}
	 that provides 
	 \begin{align*}
	 \bar Y_t 
	 &\leq  \inf_{(p,v(\cdot)) \in \mathcal S}\E^{p,v}\Big [h(X_1)-   \int_t^1\lambda_\ast(X_s,p_s,v_s(\cdot)) ds\Big|\mathcal F_t \Big], \quad t \in [0,1]. 
	 \end{align*}
	 by the arbitrariness of $(p, v(\cdot) \in \mathcal S$.

\medskip 

\noindent $(\Rightarrow)$ Again by Proposition \ref{P:74}, we may choose, for any $n\geq 1$, $(p_t^n, v_t^n(\cdot))_{t \in [0,1]}$ such that $|p_s^n|_\Xi \leq L_z$, $|1-v_s^n|_{L^2(\nu_1)}, \leq L_u$, and 
\begin{equation*}
	\lambda(X_t,\bar Z_t,\bar U_t(\cdot)) \geq -\bar Z_t p^n_t - \langle \bar U_t, 1-v_t^n \rangle_{L^2(\nu_1)}-  \lambda_\ast(X_t,p^n_t,v^n_t(\cdot)) - \frac{1}{n}.   
\end{equation*}
By a measurable selection theorem (see e.g. Theorems 8.2.10 and 8.2.11 in \cite{AuFra})  we can assume that $(p^n_t)_{t \in [0,1]}$ are progressively measurable and 
$(v^n_t(\cdot))_{t \in [0,1]}$ are predictable. 

Recalling \eqref{BSDEpv}, we get 
	\begin{align*}
\bar Y_t &= h(X_1)+   \int_t^1\lambda(X_s,\bar Z_s,\bar U_s(\cdot)) ds 
	 - \int_t^1 \bar Z_s d  W^1_s - \int_t^1 \int_{H \setminus 0} \bar U_s(w) [N_1(ds \, dw) - \nu_1(dw) ds]\\
	 & \geq h(X_1)-\int_t^1 \bar Z_s p^n_s ds - \int_t^1 \int_{H \setminus 0}\bar U_s(w) (1-v_s^n(w)) \nu_1(dw) ds-  \int_t^1\lambda_\ast(X_s,p^n_s,v^n_s(\cdot))ds\\
	 & - \frac{1-t}{n}- \int_t^1 \bar Z_s d  W^1_s - \int_t^1 \int_{H \setminus 0} \bar U_s(w) [N_1(ds \, dw) - \nu_1(dw) ds]\\
	 & = h(X_1) -  \int_t^1\lambda_\ast(X_s,p^n_s,v^n_s(\cdot))ds- \frac{1-t}{n}\\
	 & - \int_t^1 \bar Z_s d  W^{p^n}_s - \int_t^1 \int_{H\setminus 0} \bar U_s(w) [N_1(ds \, dw) - v_s^n(w)\nu_1(dw) ds]. 
\end{align*}
Applying the conditional expectation with respect to $\P^{p^n,v^n}$ we get 
	\begin{align*}
\bar Y_t + \frac{1-t}{n}
	 &\geq  \E^{p^n,v^n}\Big [h(X_1)-   \int_t^1\lambda_\ast(X_s,p^n_s,v^n_s(\cdot)) ds\Big|\mathcal F_t \Big], \quad t \in [0,1]. 
	 \end{align*}
We can therefore conclude that 
	\begin{align*}
 \bar Y_t 
	 &\geq  \liminf_{n \rightarrow \infty}\E^{p^n,v^n}\Big [h(X_1)-   \int_t^1\lambda_\ast(X_s,p^n_s,v^n_s(\cdot)) ds\Big|\mathcal F_t \Big]\\
	 &\geq  \inf_{(p,v(\cdot)) \in \mathcal S} \E^{p,v}\Big [h(X_1)-   \int_t^1\lambda_\ast(X_s,p_s,v_s(\cdot)) ds\Big|\mathcal F_t \Big], \quad t \in [0,1]. 
	 \end{align*}
\endproof

\renewcommand\thesection{\sc Appendix}
\section{}
\renewcommand\thesection{\Alph{subsection}}
\renewcommand\thesubsection{\Alph{subsection}}

\subsection{Girsanov Theorems}\label{AppGirs}
In the sequel we will denote by $\mathcal E(S)$ the Doléans-Dade exponential of  general real-valued semimartingale $S$. We recall that 
\begin{equation}\label{DDexp}
\mathcal E_T\left(S\right)= e^{S_T-\frac{1}{2} \langle S^c\rangle_T }\prod_{s \leq T}(1+ \Delta S_s)e^{-\Delta S_s}.
\end{equation}
We recall the following version of the Girsanov theorem for martingales, see e.g. Theorem 15.3.10 in \cite{cohen:stochcalculus}.

\begin{theorem}\label{Girsanov_gen}
Let $W$ be an $(\mathcal F_t)$ {cylindrical Wiener process} on $\Xi$,  and $\mu(\omega, dt\,d\ell)$ a locally integrable integer-valued  $(\mathcal F_t)$-random measure on  $\Xi$, with compensator    $\nu(\omega, dt\,d\ell)$.  Assume that $\nu(\{t\} \times \Xi)=0$ up to indistinguishability (so $\mu$ has no predictable jumps). %such that 
%$$
%\E\Big[\int_{\Xi \setminus 0}||\ell||_{\Xi}^2 \Lambda(t, \ell) F(d\ell)\Big] < \infty, \quad \forall t \in [0,T]. 
%$$

Let $\beta: \Omega \times [0,\,T] \rightarrow \Xi^\ast$ and  $\Gamma: \Omega \times [0,\,T] \times \Xi\rightarrow \R_+$ be predictable functions. Suppose that $\Gamma-1$ is stochastically $\nu$-integrable, and set 
%such that $\Gamma(\omega, t, \cdot) -1\in L^2(F)$ for all $(\omega, t) \in \Omega \times [0,\,T]$. Set  
\begin{equation}\label{MGir}
	M:=\int_{[0,\,\cdot]}\beta( t) d W_t + \int_{]0,\,\cdot]} \int_{\Xi \setminus 0}[\Gamma( t, \ell)-1] (\mu( dt\,d\ell)-\nu(dt\,d\ell)). 
\end{equation} 
	If $\mathcal E(M)$ is a uniformly integrable martingale,   then under  $\tilde \P$ defined by 
	$$
	\frac{d \tilde \P}{d\P}:= \mathcal E_T\left(M\right)
	$$
	the process
	$
	W - \int_{[0,\,\cdot]}\beta( t) dt
	$
	is a cylindrical Wiener process, and the compensator of $\mu(\omega, dt\,d\ell)$ is $\Gamma(\omega, t, \ell)\nu(\omega, dt\,d\ell)$.
\end{theorem}
\begin{remark}
Noting that  by  \eqref{MGir}
%$$
%1+ \Delta M_s =  \int_{\Xi \setminus 0}\Gamma(s, \ell) N(\{s\}\times d\ell)
%$$
%and 
$$
-\Delta M_s =  \int_{\Xi \setminus 0}(1-\Gamma(s, \ell)) N(\{s\}\times d\ell), 
$$
we have 
$$
 \prod_{s \leq T}e^{-\Delta M_s}= e^{\sum_{s \leq T}(-\Delta M_s)}= e^{\int_{]0,T]}\int_{\Xi \setminus 0}(1-\Gamma(s, \ell)) \mu(ds\, d\ell)},
$$
so \eqref{DDexp} reads 
	\begin{align*}
		\mathcal E_T\left(M\right)&=e^{\int_0^T\beta(t) d W_t-\frac{1}{2} \int_0^T|\beta(t)|^2 dt}e^{\int_{]0,\,T]} \int_{\Xi \setminus 0}[1-\Gamma(t, \ell)] \nu(d\ell\,dt)}\prod_{n \geq 1: T_n \leq T}\Gamma(T_n, \xi_{n}), 
	\end{align*}
where $(T_n, \xi_{n})_n$ denotes the sequence of jump times and positions associated to the counting measure $\mu$.  
\end{remark}
The following result provides a sufficient condition under which a local martingale $M$ is uniformly integrable, see e.g. Corollary 15.4.4 in \cite{cohen:stochcalculus}. 
\begin{corollary}\label{C_1543}
	Let $M$ be a local martingale with $\Delta M \geq -1$, and let 
	$$
	T = \inf\{t: \Delta M_t =-1\}= \inf\{t: \mathcal E_t(M) =0\}. 
	$$
	If $M$ is a locally square integrable local martingale, and 
	$$
	\E\Big[e^{\frac{1}{2} \langle M^c\rangle_T +   \langle M^d\rangle_T}\Big]< \infty,  
	$$
	then $\mathcal E(M)$ is a uniformly integrable martingale and $\{\mathcal E_\infty(M) >0 \}= \{T= \infty\}$ almost surely. 
\end{corollary}

With the help of Corollary \ref{C_1543}, when $\mu$ is a Poisson random measure one can specialize   Theorem \ref{Girsanov_gen} in the following way. 
\begin{theorem}\label{Girsanov}
Let $(\Omega, \mathcal F, (\mathcal F_t),\P)$ be a fixed filtered probability space. 
Let $W$ be an $(\mathcal F_t)$ cylindrical Wiener process on $\Xi$,  and $N(dt\,d\ell)$ an $(\mathcal F_t)$-Poisson random measure on  $\Xi$, with compensator $F(d\ell)ds$, such that $\int_{\Xi \setminus 0}||\ell||_{\Xi}^2 F(d\ell) < \infty$. 
	Suppose we have uniformly bounded functions
	 $\beta: \Omega \times [0,\,T] \rightarrow \Xi^\ast$, $\Gamma: \Omega \times [0,\,T] \times \Xi\rightarrow \R_+$ such that $\Gamma(\omega, t, \cdot) -1\in L^2(F)$ for all $(\omega, t) \in \Omega \times [0,\,T]$, and define 
		$$
	M:=\int_{[0,\,\cdot]}\beta(t) d W_t + \int_{]0,\,\cdot]} \int_{\Xi \setminus 0}[\Gamma(t, \ell)-1] (N(dt\,d\ell)-F(d\ell)dt). 
	$$ 
	Then, setting 
	$$
	\frac{d \tilde \P}{d\P}:= \mathcal E_T\left(M\right), 
	$$
	we have that  $\frac{d \tilde \P}{d\P}\Big|_{\mathcal F_t}$ is a positive square integrable martingale. Moreover, under $\tilde \P$, 
	$$
	W^{\tilde \P}:= W - \int_{[0,\,\cdot]}\beta(t) dt
	$$
	is a cylindrical Wiener process, and the compensator of $N(\omega, dt\,d\ell)$ is $\Gamma(\omega, t, \ell)F(d\ell) dt$.
\end{theorem}

%\subsection{Tanaka formula} We recall the Tanaka formula for semimartingale, see e.g. \cite{cohen:stochcalculus}.
%\begin{lemma}
%Let $(\Omega, \mathcal F,\P)$ be a fixed probability space. 
%	Let $X$ be a semimartingale, and $a \in \R$. Then there exists a continuous increasing local time $L^a$ and a pure jump process $J^{X,a}$ (unique $\P$-a.s.), such that  $J^{X,a}(0)=0$ and 
%	$$
%	d |X_t - a|\leq \textup{sign}(X_{t-}-a) dX_t + dL_t^a + \Delta J^{X,a}
%	$$
%	where 
%	$$
%	 \Delta J^{X,a} := |X_t -a |-|X_{t-}-a|-\textup{sign}(X_{t-}-a)\Delta X_t \geq 0.
%	$$
%\end{lemma}

\bibliographystyle{plainnat}
\bibliography{Bibliography}
\end{document}